\documentclass[11pt]{article}
\usepackage{amsfonts}
\usepackage{amssymb}
\usepackage{amsthm}
\usepackage{amsmath}
\usepackage{graphicx}
\usepackage{empheq}
\usepackage{indentfirst}
\usepackage{cite}
\usepackage{mathrsfs}
\usepackage{graphics}

 \newtheorem{theorem}{Theorem}[section]
 
 \newtheorem{lemma}{Lemma}[section]
 \newtheorem{proposition}{Proposition}[section]

 \newtheorem{remark}{Remark}[section]
 \numberwithin{equation}{section}

\allowdisplaybreaks[4]

\def\ddiv{\mathrm{div}}

\def\e{\epsilon}

\newcommand{\beq}{\begin{equation}}
\newcommand{\eeq}{\end{equation}}
 \def\non{\nonumber }
\def\bea{\begin{eqnarray}}
\def\eea{\end{eqnarray}}
\def\tL{\textbf{L}^2(\mathbb{R}^3)}

\def\tu{\tilde{u}}
\def\tv{\tilde{v}}
\def\tpsi{\tilde{\psi}}

\begin{document}

\title{Global solution to the drift-diffusion-Poisson system for semiconductors with nonlinear recombination-generation rate }
\author{Hao Wu\thanks{School of Mathematical Sciences and Shanghai Key Laboratory for Contemporary Applied
Mathematics,
Fudan University, Shanghai 200433,\ P.R. China,
\textsl{haowufd@yahoo.com}.} \ and \
Jie Jiang\thanks{Wuhan Institute of Physics and Mathematics, Chinese Academy of Sciences, Wuhan 430071, HuBei Province, P.R. China,
\textsl{jiangbryan@gmail.com}.}
}

\date{\today}
\maketitle
\begin{abstract}
In this paper, we study the Cauchy problem of a time-dependent drift-diffusion-Poisson system
for semiconductors. Existence and uniqueness of global weak solutions are proven for the system with a higher-order nonlinear recombination-generation
rate $R$. We also show that the global weak solution will converge to a unique equilibrium as time tends to infinity. \\

\noindent  {\bf Keywords}: drift-diffusion-Poisson system; global weak solution; uniqueness; long-time behavior.\\
\textbf{AMS Subject Classification}: 35G25, 35J20, 35B40, 35B45
\end{abstract}

\section{Introduction}\setcounter{equation}{0}
We consider the following drift-diffusion-Poisson model for semiconductors that is a coupled system of parabolic-elliptic equations:
 \beq
 \label{DDP}
 \begin{cases}
 n_t=\ddiv\big(\nabla n+n\nabla(\psi+V_n)\big)-R(n,p,x),\\
 p_t=\ddiv\big(\nabla p+p\nabla(-\psi+V_p)\big)-R(n,p,x),\\
 -\varepsilon^2\Delta\psi=n-p-D(x).
 \end{cases}
 \eeq
 System \eqref{DDP} models the transport of the electrons and holes in semiconductor and plasma devices (cf. \cite{MRS,J01}).
 $n=n(x,t)$ is the spatial distribution of electrons (negatively charged) and $p=p(x,t)$ is the spatial
distribution of holes (positively charged). $\psi=\psi(x,t)$ is
the self-consistent electrostatic potential created by the two
charge carrier species (electrons and holes) and by the doping
profile $D=D(x)$ of the semiconductor device. The
charge carriers are assumed to be confined by the external potentials $V_n$ and $V_p$. This
replaces the usual assumption of a bounded domain (cf.
\cite{FI,Ga,MRS} and the references
therein). The function $R = R(n, p, x)$ represents the so-called recombination-generation rate for electrons and holes.
 The parameter $\varepsilon$ appearing in the Poisson
equation is the scaled Debye length of the semiconductor device that stands for the screening of the hole and electron particles.
In this paper, we are interested in the Cauchy problem to system \eqref{DDP} and assume that \eqref{DDP} is subject to
 the following initial data
 \beq
\label{initial}
 n(x,t)|_{t=0}=n_I(x)\geq 0,\qquad p(x,t)|_{t=0}=p_I(x)\geq 0.
 \eeq

The generation and recombination of electrons and holes in a semiconductor play an important role
in their electrical and optical behavior \cite{D}. Recombination is a process by which both carriers annihilate
each other: the electrons fall in one or multiple steps into the empty state that is
associated with the hole. Generation can be viewed as its inverse process whereby electrons and holes are created. There are several typical recombination mechanisms that the energy of carriers will be dissipated during these processes by different ways (cf. e.g., \cite{MRS, J01,D}). For instance,
\begin{enumerate}
\item {\emph{Band-to-band recombination} (also referred to as direct thermal recombination). The energy is emitted in the form of a photon.
The recombination rate depends on the density of available electrons and holes and it can be expressed as
\beq
R(n,p)=C(np-n_i^2),\label{bb}
\eeq
where $n_i$ denotes the intrinsic carrier density of the
semiconductor. }

\item {\emph{Shockley-Read-Hall (SRH) recombination} (also called the trap-assisted recombination). A
two-step transition of an electron from the conduction band to the valence
band occurs and $R$ is in the form
\beq
R(n,p)=\frac{(np-n_i^2)}{r_1n+r_2p+r_3},\label{rhs}
\eeq
where $r_1,r_2,r_3$ are proper positive functions. }

\item {\emph{Auger recombination}. An electron and a hole recombine in a
band-to-band transition, but the resulting energy is given off to another electron
or hole in the form of kinetic energy. The corresponding recombination rate is similar to that of band-to-band
recombination, but involves a third particle:
\beq
R(n,p)=(C_n n+C_p p)(np-n_i^2).\label{aug}
\eeq}
\end{enumerate}

Extensive mathematical study of the drift-diffusion-Poisson system has been developed in the literature.
 For the initial boundary value problem of \eqref{DDP} in a bounded domain $\Omega\subset \mathbb{R}^N$ with
  various boundary conditions (e.g., the Neumann type, Robin type, or mixed boundary conditions), existence and uniqueness
    as well as long-time behavior have been investigated by many authors, see for instance,
    \cite{AMV,BD,BHN,FI1,FI,FI2,GGH94, GGH, GH97, GH2005, Glitzky08, Glitzky09, Glitzky, Gro87,GW91,S,XW1} and reference therein.
    On the other hand, for the sake of modeling simplicity and for the particularly
interesting mathematical features, it would also be interesting to
consider the Cauchy problem of \eqref{DDP} (cf. e.g., \cite{AMT, BDM, FFM, KK, KO}). Existence and uniqueness
results and stability of strong solutions in $L^p(\mathbb{R}^N)$
spaces $(N\geq 2)$ were proven in \cite{KO} for a system analogous
to \eqref{DDP}. However, in their system there were no external
potentials and the recombination-generation rate $R$ was replaced by
a given function $f=f(x,t)\in L^{\theta}(0,T;W^{1,\eta})$ with $1\leq
\theta<2$, $\frac{N}{2}<\eta<N$, which expressed the variation of the
charge by the external current.
 As far as the long-time behavior of global solutions to the Cauchy problem is concerned, when the recombination-generation
 term $R$ is absent, exponential convergence to equilibrium with a confining potential and an algebraic rate towards a self-similar state without confinement have been obtained in \cite{AMT}. The analysis therein is based on the well-known entropy approach for diffusion and diffusion-convection equations that has been extensively studied in recent years (cf. \cite{AMTU,CJMTU,UAMT} and the references therein). We also refer to \cite{KK} in which an optimal $L^p$ decay estimate of solutions to the Cauchy problem was obtained via a time weighted energy method (without confinement and recombination-generation rate $R$). When the recombination-generation process is taken into account, the situation is more complicated.
 In \cite{WMZ}, the authors proved the global existence and uniqueness of weak solutions of problem \eqref{DDP} in $\mathbb{R}^3$ with an (unbounded) external confining potential $V_n=V_p=V$ and under the restrictive assumption that $R$ has a linear growth (which, however, recovers the Shockley-Read-Hall recombination, cf. \eqref{rhs}). Besides, existence and uniqueness of the steady state and partial result on the convergence to equilibrium were obtained. Recently, exponential $L^1$ convergence to equilibrium was proved in \cite{FFM} via entropy method  for global solutions to a simplified convection-diffusion-reaction model with confinement and Shockley-Read-Hall recombination-generation rate but neglecting the influence of the self-consistent potential $\psi$. It would be interesting to study the well-posedness as well as long time behavior of the full convection-diffusion-reaction-Poisson system \eqref{DDP}--\eqref{initial} with more general recombination-generation rate $R$ including the higher nonlinear cases \eqref{bb} and \eqref{aug}.

For the sake of simplicity, we consider the whole-space case posed
on $\mathbb{R}^3$. Similar results can be obtained for the
two-dimensional whole space case with some minor modifications, due
to the different properties of the Newtonian potential.

 In order to
formulate our assumptions and results, we first introduce some notations on the functional settings.
$H^m(\mathbb{R}^3)$ ($m\in\mathbb{N}$) is
 used to denote the Sobolev
space $W^{m,2}(\mathbb{R}^3)$, and $\|\cdot\|_{H^m(\mathbb{R}^3)}$
is its corresponding norm. We denote $L^r(\mathbb{R}^3)$ $(r\geq 1)$
with norm $\|\cdot\|_{L^r(\mathbb{R}^3)}$ and the vector space
$\mathbf{L}^r(\mathbb{R}^3)= (L^r(\mathbb{R}^3))^3$ $(r\geq 1)$ with
norm $\|\cdot\|_{\mathbf{L}^r(\mathbb{R}^3)}$. Moreover, for a
potential function $V$, we define the weighted $L^r$ ($r>1$) space
as follows \beq L^r(\mathbb{R}^3, e^{V(x)}dx):=\left\{u\in
L^1_{\text{loc}}(\mathbb{R}^3)\left|\;\int_{\mathbb{R}^3}|u(x)|^r
e^{V(x)}dx<\infty\right.\right\} \eeq with norm
$\|u\|_{L^r(\mathbb{R}^3,
e^{V(x)}dx)}:=\left(\int_{\mathbb{R}^3}|u(x)|^r
e^{V(x)}dx\right)^{1/r}$. We define the weighted vector space and
 its norm  in a similar way, which are denoted by $\mathbf{L}^2(\mathbb{R}^3, e^{V(x)}dx)$ and $\|\cdot\|_{\mathbf{L}^2(\mathbb{R}^3, e^{V(x)}dx)}$, respectively. For any Hilbert space
$H$, we denote its subspace
$$ H_+=\{f(x)\in H\ |\  f(x)\geq 0, \ \text{a.e.}\  x\in \mathbb{R}^3\}.$$
Throughout this paper, we use $C,\, C_i (i\in\mathbb{N})$ to
 denote genetic constants that may vary in different places (even in the same
 estimate). Particular dependence of those constants will be
 explained in the text if necessary.

 Next, we make the following assumptions on confining potentials $V_n, V_p$, the recombination-generation rate $R$ and the doping profile $D$:
\begin{itemize}
\item[(H1a)] There exist constants $\rho_n, \rho_p>0$ such that
 \beq
 \frac{\partial^2 V_n}{\partial x^2}\geq \rho_n \mathbb{I},\ \
 \frac{\partial^2 V_p}{\partial x^2}\geq \rho_p \mathbb{I},\quad \forall x\in \mathbb{R}^3,\non
 \eeq
in the sense of positive-defined matrix. Moreover, there exists $K>0$ such that
 \beq
 \label{V1}
 \|\Delta V_n\|_{L^{\infty}(\mathbb{R}^3)}\leq K, \quad \|\Delta V_p\|_{L^{\infty}(\mathbb{R}^3)}\leq K.
 \eeq
 \item[(H1b)] There exists $K'>0$ such that
  \beq
   \|V_n(x)-V_p(x)\|_{L^\infty(\mathbb{R}^3)}\leq K'<+\infty.\label{V2}
  \eeq
 \item[(H2a)] The recombination-generation rate $R=R(n,p,x)$ is of the form
      \beq
      R(n,p,x)=F(n,p)\left(np-\delta^2\mu_n\mu_p\right),\non
      \eeq
      where $\mu_n(x)=e^{-V_n(x)}$, $\mu_p(x)=e^{-V_p(x)}$. $\delta$ is a positive constant standing for the scaled average intrinsic carrier density of the
      semiconductor.
      Without loss of generality, we assume that $\delta=1$ in the remaining part of this paper.
\item[(H2b)] The scalar function $F:\mathbb{R}^2\rightarrow \mathbb{R}$ is a Lipschitz continuous function with linear growth, namely, there exist constants $c_1, c_2>0$
 independent of $n, p$ such that
 \beq
 |F(n_1,p_1)-F(n_2,p_2)|\leq c_1(|n_1-n_2|+|p_1-p_2|),\ \ \forall \, n_1,p_1,n_2,p_2 \in \mathbb{R},\non
 \eeq
 \beq
 |F(n,p)|\leq c_2(1+|n|+|p|)\quad\forall \ n,p\in \mathbb{R}.\non
 \eeq
Moreover, $F(n,p)\geq 0$ if $n, p\geq 0$.

\item[(H3)] $ D(x)\in L^1(\mathbb{R}^3)\cap L^{\infty}(\mathbb{R}^3).$
\end{itemize}

\begin{remark}\label{Rem1}
It easily follows from (H1a) that the confining potentials $V_n(x)$ and $V_p(x)$ are uniformly convex and can be bounded from below by a finite number $V_b\in \mathbb{R}$.
 Thus, $\|\mu_n\|_{L^\infty}\leq e^{-V_b}$ and $\|\mu_p\|_{L^\infty}\leq e^{-V_b}$.
  Without loss of generality, we assume in \eqref{DDP} that the diffusion coefficients and carrier mobilities are equal to one. Moreover, for the sake of simplicity, we set
 \beq
 \int_{\mathbb{R}^3}\mu_n dx=\int_{\mathbb{R}^3}\mu_p dx=1.\non
 \eeq
 The above simplifications do not affect the subsequent analysis.
 We also infer from \eqref{V2} that the norms on $L^2\left(\mathbb{R}^3,
 e^{V_n(x)}dx\right)$ and $L^2\left(\mathbb{R}^3,
 e^{V_p(x)}dx\right)$ are equivalent.

 A typical example for the confining potential is $\frac{|x|^2}{2}$ (cf. e.g., \cite{AMT,AMTU, FFM}). We remark that the confining potential is introduced in the Cauchy problem (although somewhat unphysical) in order to prevent the particles from escaping to infinity as time progresses. For the Cauchy problem of drift-diffusion-Poisson system without generation-recombination term,
 different types of large time behavior for the cases with or without confining potential have been illustrated in \cite{AMT} (see also \cite{BDM} for a related system modeling the bipolar plasma).
 \end{remark}

Now we state the the main results of this paper.

\begin{theorem}[Well-posedness] \label{Global} Suppose that (H1a)--(H3) are satisfied.
 For any initial data $n_I\in L^2(\mathbb{R}^3, e^{V_n(x)}dx)\cap L^\infty(\mathbb{R}^3)$, $p_I\in L^2(\mathbb{R}^3, e^{V_p(x)}dx)\cap L^\infty(\mathbb{R}^3)$, $n_I,p_I\geq 0$, problem \eqref{DDP}--\eqref{initial} admits a unique global weak solution $(n,p,\psi)$ such that for any $T>0$,
 \bea
 && n\in L^\infty(0,T; L^2(\mathbb{R}^3, e^{V_n(x)}dx)\cap L^\infty(\mathbb{R}^3)),\ \nabla n\in L^2(0,T;\mathbf{L}^2(\mathbb{R}^3, e^{V_n(x)}dx)),\non
 \\
 && p\in L^\infty(0,T; L^2(\mathbb{R}^3, e^{V_p(x)}dx)\cap L^\infty(\mathbb{R}^3)),\ \nabla p\in L^2(0,T;\mathbf{L}^2(\mathbb{R}^3, e^{V_p(x)}dx)),\non
 \\
 && n_t,p_t\in L^2(0,T; (H^1(\mathbb{R}^3))'),\non\\
 && n(t)\geq0, \quad p(t)\geq0,\quad t\in[0,T],\quad \text{ a.e. } x\in\mathbb{R}^3,\non\\
 &&\nabla \psi\in L^\infty(0,T;\mathbf{L}^\infty(\mathbb{R}^3)),\  \Delta \psi\in L^\infty(0,T;L^\infty(\mathbb{R}^3)),\non
 \eea
 where $\psi=\psi(x,t)$ is given by the Newtonian potential
 \beq
 \psi(x,t)=\frac{1}{S_3}\int_{\mathbb{R}^3}\frac{n(y,t)-p(y,t)-D(y)}{|x-y|}dy,\non
 \eeq
 with $S_3=\frac{2\pi^\frac32}{\Gamma(\frac32)}$ being the surface area of the 2D unit ball.
\end{theorem}

\begin{theorem}[Long-time behavior]
 \label{connn}
 Under the assumptions in Theorem \ref{Global}, we have the following uniform-in-time estimate for the global solutions:
 \beq
  \|n(t)\|_{L^\infty(0,+\infty; L^r(\mathbb{R}^3))}+\|p(t)\|_{L^\infty(0,+\infty; L^r(\mathbb{R}^3))}<\infty,\quad\forall\, r\in[1,+\infty]. \label{unie}
  \eeq
 Moreover, for every fixed $t^*>0$, the global shifted solution $(n(t+s), p(t+s), \psi(t+s))$ $(s\in (0,t^*) )$ of problem \eqref{DDP}--\eqref{initial}
 converges to the unique steady state $(n_\infty,p_\infty,\psi_\infty)$ that satisfies \eqref{sta} as $t\to +\infty$ in the
 following sense:
 \bea && n(t+\cdot)\rightarrow n_\infty,\quad p(t+\cdot)\rightarrow p_\infty\quad  \text{in} \
 L^1((0,t^*)\times\mathbb{R}^3),
 \non\\
 && \nabla \psi(t+\cdot) \rightharpoonup \nabla \psi_\infty\quad
\text{in} \ L^2(0,t^*; \mathbf{H}^1(\mathbb{R}^3)),
 \non\\
 && \psi(t+\cdot) \rightharpoonup \psi_\infty\quad  \text{in} \
L^2(0,t^*; L^6(\mathbb{R}^3)).\non
 \eea
 \end{theorem}

\begin{remark}
 Our results holds for arbitrary but fixed $\varepsilon>0$. In the following analysis, we
just set $\varepsilon=1$ without of loss of generality.
  We  note that the quasi-neutral limit (namely, zero-Debye-length limit
 $\varepsilon\to 0$) of the drift-diffusion-Poisson system is a challenging and physically
complex modeling problem for bipolar kinetic models of semiconductors, which has
been analyzed by several authors, see, e.g. \cite{HW, WXM} and the references
cited therein.
\end{remark}

As we have mentioned before, for a class of recombination-generation
rate with at most linear growth, existence of global weak solutions
to  problem \eqref{DDP}--\eqref{initial} in $\mathbb{R}^3$ has been
obtained in \cite{WMZ}. However, argument therein fails to apply in
our present case due to the higher-order nonlinear reaction term $R$
that includes both the band-to-band and the Auger recombination (cf.
(H2a)--(H2b)). On the other hand, well-posedness results for
drift-diffusion-Poisson system with higher-order
recombination-generation rate have been proved in the bounded domain
in $\mathbb{R}^N$, $N\leq3$ (see \cite{FI,FI1,FI2} for the case with
$F$ being bounded, and \cite{XW1} for the case that $F$ has a linear
growth). Since now we are considering the Cauchy problem in the
whole space that is unbounded and the carriers are confined by
unbounded external potentials, the methods for the initial boundary
value problem cannot be applied directly. We need to exploit and
employ several techniques in the literature to prove the global
existence and uniqueness of solutions to problem
\eqref{DDP}--\eqref{initial}.
 In order to overcome the difficulties from the higher-order reaction term $R$, we first introduce a $L^\infty$ cut-off to the unknowns $n,p$ in $R$ and study an approximation problem associated with our original system \eqref{DDP}--\eqref{initial}. To deal with the unbounded confining potentials, we then transform
the approximate problem into a new form by introducing some new variables with proper weight functions.
After obtaining the well-posedness of the approximate problem, we try to derive proper uniform estimates based on a Stampacchia-type $L^\infty$ estimation technique (cf. \cite{FI}) that enable us to pass to limit and show the existence of global weak solutions to the original system \eqref{DDP}--\eqref{initial}.
Finally, we get uniform-in-time $L^r$ ($r\in[1,+\infty]$)
estimates for the global solutions under more general assumptions by extending the methods in \cite{WMZ,FFM} and investigate the long-time behavior of global solutions.

The remaining part of the paper is organized as follows. In Section 2, we prove
the well-posedness of an approximate problem and obtain some uniform estimates that are independent of the approximate parameter.
In Section 3, we prove the existence of global solutions to the original problem \eqref{DDP}--\eqref{initial} by passing to the limit and show the uniqueness of the solution.
 In Section 4, we obtain some uniform-in-time estimates of the solutions and show that as time tends to infinity the global solutions will
converge to a unique steady state.

\section{Well-posedness of the Approximate System}\setcounter{equation}{0}
In order to overcome the difficulty brought by the higher-order
nonlinearity $R$, we introduce and study the following approximate
problem in this section. For any $\sigma>0$, consider
 \beq
 \label{APX}
 (AP)\left\{\begin{array}{ll}\partial_t n_{\sigma}=\ddiv\big(\nabla n_{\sigma}+n_{\sigma}\nabla(\psi_{\sigma}+V_n)\big)- \tilde{R}(n_\sigma, p_\sigma,x),\\
 \partial_t p_\sigma=\ddiv\big(\nabla p_\sigma+p_\sigma\nabla(-\psi_\sigma+V_p)\big)- \tilde{R}(n_\sigma, p_\sigma,x),\\
 -\Delta\psi_\sigma=n_\sigma-p_\sigma-D(x),\end{array}\right.
  \eeq
  subject to the initial data
\beq
\label{initial2}
 n_\sigma(x,t)|_{t=0}=n_I,\qquad p_\sigma(x,t)|_{t=0}=p_I.
\eeq
The approximated recombination-generation rate $\tilde{R}$ in \eqref{APX} is given by
 \bea
  \tilde{R}(n_\sigma, p_\sigma,x)&=& R\left(\frac{n_\sigma}{1+\sigma n_\sigma},\frac{p_\sigma}{1+\sigma p_\sigma},x\right)\non\\
  &=& F\left(\frac{n_\sigma}{1+\sigma n_\sigma},\frac{p_\sigma}{1+\sigma p_\sigma}\right)\left(\frac{n_\sigma}{1+\sigma n_\sigma}\frac{p_\sigma}{1+\sigma p_\sigma}-\mu_n\mu_p\right).
 \eea
 Now we state the main result of this section.
 \begin{theorem} \label{apx}
Suppose that assumptions (H1a)--(H3) are satisfied. For any $\sigma>0$, $n_I\in L^2(\mathbb{R}^3, e^{V_n(x)}dx)\cap L^4(\mathbb{R}^3)$, $p_I\in L^2(\mathbb{R}^3, e^{V_p(x)}dx)\cap L^4(\mathbb{R}^3)$, $n_I,p_I\geq 0$,
 problem \eqref{APX}--\eqref{initial2} admits a unique global weak solution $(n_\sigma,p_\sigma,\psi_\sigma)$ such that for any $T>0$,
 \bea
 && n_\sigma\in C([0,T]; L^2(\mathbb{R}^3, e^{V_n(x)}dx)),\quad \nabla n_\sigma\in L^2((0,T);\mathbf{L}^2(\mathbb{R}^3, e^{V_n(x)}dx)),\non\\
 && p_\sigma\in C([0,T]; L^2(\mathbb{R}^3, e^{V_p(x)}dx)),\quad \nabla p_\sigma\in L^2((0,T);\mathbf{L}^2(\mathbb{R}^3, e^{V_p(x)}dx)),\non
 \eea
 \beq
 n_\sigma(t)\geq 0, \quad p_\sigma(t)\geq 0,\quad t\in[0,T],\quad \text{a.e.} \ \ x\in\mathbb{R}^3.\non
 \eeq
 $\psi_\sigma=\psi_\sigma(x,t)$ is the Newtonian potential with respect to $n_\sigma(x,t), p_\sigma(x,t)$ given by
 \beq
 \psi_\sigma(x,t)=\frac{1}{S_3}\int_{\mathbb{R}^3}\frac{n_\sigma(y,t)-p_\sigma(y,t)-D(y)}{|x-y|}dy,\quad\text{with}\ S_3=\frac{2\pi^\frac32}{\Gamma(\frac32)}.\non
 \eeq
 \end{theorem}
The proof of Theorem \ref{apx} consists of several steps. First, we derive some properties for the new reaction term $\tilde{R}$ under assumptions (H2a)--(H2b).

\begin{lemma}\label{lipc}
Under assumptions (H2a)--(H2b),
the function $\tilde{R}=\tilde{R}(n_\sigma,p_\sigma,x)$
satisfies the following properties \\
 (i) $\tilde{R}$ has at most a linear growth for any $n_\sigma, p_\sigma\geq 0$, i.e.,
 \beq
 |\tilde{R}(n_\sigma,p_\sigma,x)|\leq C_\sigma(a(x)+n_\sigma+p_\sigma),
 \eeq
 where $0\leq a(x) \in L^1(\mathbb{R}^3)\cap \ L^\infty
 (\mathbb{R}^3)$ and $a(x)\in L^2_+\left(\mathbb{R}^3, e^{V_n(x)}dx\right)\cap L^2_+\left(\mathbb{R}^3, e^{V_p(x)}dx\right)$.\\
 (ii) $\tilde{R}$ is Lip-continuous in $L^2_+\left(\mathbb{R}^3, e^{V_i(x)}dx\right)$, $i=\{n,p\}$,
 such that
 for any $n_{\sigma}^{(1)}$, $n_{\sigma}^{(2)}$  $\in L^2_+\left(\mathbb{R}^3,
 e^{V_n(x)}dx\right)$,  $p_{\sigma}^{(1)}$, $p_{\sigma}^{(2)}$ $\in L^2_+\left(\mathbb{R}^3,
 e^{V_p(x)}dx\right)$, it holds
 \bea  && \|\tilde{R}(n_{\sigma}^{(1)},p_{\sigma}^{(1)},\cdot)-\tilde{R}(n_{\sigma}^{(2)},p_{\sigma}^{(2)},\cdot)\|_{L^2(\mathbb{R}^3,
 e^{V_i(x)}dx)}\non\\
 &\leq& \tilde{K}
 \left(\|n_{\sigma}^{(1)}-n_{\sigma}^{(2)}\|_{L^2(\mathbb{R}^3, e^{V_i(x)}dx)}+\|p_{\sigma}^{(1)}-p_{\sigma}^{(2)}\|_{L^2(\mathbb{R}^3, e^{V_i(x)}dx)}\right),\label{LIP}
 \eea
 where the constant $\tilde{K}$ may depend on $c_1, c_2, \sigma$ and $V_b$ in Remark \ref{Rem1}.
\end{lemma}
\begin{proof}
We observe that for any $\varphi\geq 0$, it holds
 \beq
 0\leq \frac{\varphi}{1+\sigma \varphi}\leq\frac{1}{\sigma} \quad \text{and} \ \ \frac{\varphi}{1+\sigma \varphi}\leq \varphi,\quad  \forall\,  \sigma>0.
 \label{bound}
 \eeq
Due to the above simple facts and assumptions (H2a)--(H2b), we can verify
that
 \bea
 |\tilde{R}(n_\sigma, p_\sigma, x)|&\leq& c_2\left(1+ \frac{n_\sigma}{1+\sigma n_\sigma}+\frac{p_\sigma}{1+\sigma p_\sigma}\right)\left|\frac{n_\sigma}{1+\sigma n_\sigma}\frac{p_\sigma}{1+\sigma p_\sigma}-\mu_n\mu_p\right|\non\\
 &\leq& c_2\left(1+\frac{2}{\sigma}\right)\mu_n\mu_p+ c_2\left(\frac{1}{\sigma}+\frac{1}{\sigma^2}\right)(n_\sigma+p_\sigma).\label{glin}
 \eea
 Then we can simply set $C_\sigma=c_2\left(1+\frac{2}{\sigma}+\frac{1}{\sigma^2}\right)$ and $a(x)=\mu_n\mu_p$, which obviously satisfies the required conditions by assumption (H1a).

  For any $n_{\sigma}^{(1)}$, $n_{\sigma}^{(2)}$  $\in L^2_+\left(\mathbb{R}^3,
 e^{V_n(x)}dx\right)$,  $p_{\sigma}^{(1)}$, $p_{\sigma}^{(2)}$ $\in L^2_+\left(\mathbb{R}^3,
 e^{V_p(x)}dx\right)$, we infer from \eqref{bound} that
 \beq
 \left|\frac{n_{\sigma}^{(1)}}{1+\sigma n_{\sigma}^{(1)}}-\frac{n_{\sigma}^{(2)}}{1+\sigma n_{\sigma}^{(2)}} \right|\leq |n_{\sigma}^{(1)}-n_{\sigma}^{(2)}|, \quad \left| \frac{p_{\sigma}^{(1)}}{1+\sigma p_{\sigma}^{(1)}}-\frac{p_{\sigma}^{(2)}}{1+\sigma p_{\sigma}^{(2)}}\right|\leq |p_{\sigma}^{(1)}-p_{\sigma}^{(2)}|,\label{li1}
 \eeq
 and as a result,
 \bea
 && \left|\frac{n_{\sigma}^{(1)}}{1+\sigma n_{\sigma}^{(1)}} \frac{p_{\sigma}^{(1)}}{1+\sigma p_{\sigma}^{(1)}}-\frac{n_{\sigma}^{(2)}}{1+\sigma n_{\sigma}^{(2)}} \frac{p_{\sigma}^{(2)}}{1+\sigma p_{\sigma}^{(2)}}\right|\non\\
 &\leq& \frac{p_{\sigma}^{(1)}}{1+\sigma p_{\sigma}^{(1)}} \left|\frac{n_{\sigma}^{(1)}}{1+\sigma n_{\sigma}^{(1)}}-\frac{n_{\sigma}^{(2)}}{1+\sigma n_{\sigma}^{(2)}}\right|
 +\frac{n_{\sigma}^{(2)}}{1+\sigma n_{\sigma}^{(2)}} \left|\frac{p_{\sigma}^{(1)}}{1+\sigma p_{\sigma}^{(1)}}- \frac{p_{\sigma}^{(2)}}{1+\sigma p_{\sigma}^{(2)}}\right|\non\\
 &\leq& \frac{1}{\sigma}(|n_{\sigma}^{(1)}-n_{\sigma}^{(2)}|+|p_{\sigma}^{(1)}-p_{\sigma}^{(2)}|).\label{li2}
 \eea
 Denote
 \beq
  F_j=F\left(\frac{n_{\sigma}^{(j)}}{1+\sigma n_{\sigma}^{(j)}}, \frac{p_{\sigma}^{(j)}}{1+\sigma p_{\sigma}^{(j)}}\right), \quad j=1,2.\non
 \eeq
  Then we get
 \begin{eqnarray}  && \|\tilde{R}(n_{\sigma}^{(1)},p_{\sigma}^{(1)},\cdot)-\tilde{R}(n_{\sigma}^{(2)},p_{\sigma}^{(2)},\cdot)\|_{L^2(\mathbb{R}^3,
 e^{V_i(x)}dx)}\non\\
 &\leq& \left\|(F_1 -F_2) \mu_n\mu_p\right\|_{L^2(\mathbb{R}^3,
 e^{V_i(x)}dx)}\non\\
 && +  \left\|F_1 \left(\frac{n_{\sigma}^{(1)}}{1+\sigma n_{\sigma}^{(1)}} \frac{p_{\sigma}^{(1)}}{1+\sigma p_{\sigma}^{(1)}}-\frac{n_{\sigma}^{(2)}}{1+\sigma n_{\sigma}^{(2)}} \frac{p_{\sigma}^{(2)}}{1+\sigma p_{\sigma}^{(2)}}\right)\right\|_{L^2(\mathbb{R}^3,
 e^{V_i(x)}dx)}\non\\
 &&+ \left\|(F_1 -F_2) \frac{n_{\sigma}^{(2)}}{1+\sigma n_{\sigma}^{(2)}} \frac{p_{\sigma}^{(2)}}{1+\sigma p_{\sigma}^{(2)}}\right\|_{L^2(\mathbb{R}^3,
 e^{V_i(x)}dx)}\non\\
 &:=& I_1+I_2+I_3,\quad i=\{n,p\}.
 \end{eqnarray}
For the case $i=n$, it follows from \eqref{bound}, \eqref{li1}, \eqref{li2}, (H2a)--(H2b) and Remark \ref{Rem1} that
 \bea
 I_1&\leq& \left\|c_1\left(\left|\frac{n_{\sigma}^{(1)}}{1+\sigma n_{\sigma}^{(1)}}-\frac{n_{\sigma}^{(2)}}{1+\sigma n_{\sigma}^{(2)}} \right|+\left| \frac{p_{\sigma}^{(1)}}{1+\sigma p_{\sigma}^{(1)}}-\frac{p_{\sigma}^{(2)}}{1+\sigma p_{\sigma}^{(2)}}\right|\right)\mu_n\mu_p\right\|_{L^2(\mathbb{R}^3,
 e^{V_n(x)}dx)}\non\\
 &\leq& \|c_1(|n_{\sigma}^{(1)}-n_{\sigma}^{(2)}|+|p_{\sigma}^{(1)}-p_{\sigma}^{(2)}|)\mu_n\mu_p\|_{L^2(\mathbb{R}^3,
 e^{V_n(x)}dx)}\non\\
 &\leq& 2c_1e^{-2V_b} \left(\|n_{\sigma}^{(1)}-n_{\sigma}^{(2)}\|_{L^2(\mathbb{R}^3, e^{V_n(x)}dx)}+\|p_{\sigma}^{(1)}-p_{\sigma}^{(2)}\|_{L^2(\mathbb{R}^3, e^{V_n(x)}dx)}\right),\non
 \eea
 \bea
 I_2&\leq& \left\|\frac{c_2}{\sigma} \left(1+ \frac{n_{\sigma}^{(1)}}{1+\sigma n_{\sigma}^{(1)}}+ \frac{p_{\sigma}^{(1)}}{1+\sigma p_{\sigma}^{(1)}}\right)(|n_{\sigma}^{(1)}-n_{\sigma}^{(2)}|+|p_{\sigma}^{(1)}-p_{\sigma}^{(2)}|) \right\|_{L^2(\mathbb{R}^3,
 e^{V_n(x)}dx)}\non\\
 &\leq& \frac{2c_2}{\sigma}\left(1+\frac{2}{\sigma}\right)\left(\|n_{\sigma}^{(1)}-n_{\sigma}^{(2)}\|_{L^2(\mathbb{R}^3, e^{V_n(x)}dx)}+\|p_{\sigma}^{(1)}-p_{\sigma}^{(2)}\|_{L^2(\mathbb{R}^3, e^{V_n(x)}dx)}\right),\non
 \eea
 \bea
 I_3&\leq& \frac{1}{\sigma^2}\|c_1(|n_{\sigma}^{(1)}-n_{\sigma}^{(2)}|+|p_{\sigma}^{(1)}-p_{\sigma}^{(2)}|)\|_{L^2(\mathbb{R}^3,
 e^{V_n(x)}dx)}\non\\
 &\leq& \frac{2c_1}{\sigma^2} \left(\|n_{\sigma}^{(1)}-n_{\sigma}^{(2)}\|_{L^2(\mathbb{R}^3, e^{V_n(x)}dx)}+\|p_{\sigma}^{(1)}-p_{\sigma}^{(2)}\|_{L^2(\mathbb{R}^3, e^{V_n(x)}dx)}\right),\non
 \eea
 where $V_b$ is the lower bound of $V_n, V_p$ (see Remark 1.1).
 Collecting the above estimates together, we see that \eqref{LIP} $(i=n)$ holds with
 $$\tilde{K}=2c_1e^{-2V_b} +\frac{2c_2}{\sigma}\left(1+\frac{2}{\sigma}\right)+ \frac{2c_1}{\sigma^2}.$$
  The case $i=p$ can be treated in the same way. The proof is complete.
  \end{proof}
Next, for any fixed $\sigma>0,$ we introduce the following
transformation of unknown variables (cf. \cite{AMTU,WMZ}) \beq
 \label{transform}
 u(x,t):=n_\sigma(x,t) e^{\frac{V_n(x)}{2}},\qquad v(x,t):=p_\sigma(x,t)
 e^{\frac{V_p(x)}{2}},
 \eeq
 then it follows from system \eqref{APX} and a direct computation that $u$ and $v$ satisfy the following transformed approximate system
 \beq\label{TAPX}
 (TAP) \begin{cases}u_t-\Delta u+A_n(x)u=f_1(u,v,\psi_\sigma),\\
v_t-\Delta v+A_p(x) v=f_2(u,v,\psi_\sigma),\\
-\Delta\psi_\sigma= ue^{-\frac{V_n(x)}{2}}-ve^{-\frac{V_p(x)}{2}}-D(x),\end{cases}\eeq
where
 \beq
  A_n(x)=\frac{1}{4}|\nabla V_n(x)|^2-\frac{1}{2}\Delta V_n(x)+K, \ \ A_p(x)=\frac{1}{4}|\nabla V_p(x)|^2-\frac{1}{2}\Delta V_p(x)+K\non
 \eeq with $K$ being the constant in (H1a) (see \eqref{V1}).
 System  \eqref{TAPX} is
subject to the initial data
 \beq\label{initial3}
 u(x,t)|_{t=0}=n_I(x)e^{\frac{V_n(x)}{2}}:=u_I,\quad v(x,t)|_{t=0}=p_I(x)e^{\frac{V_p(x)}{2}}:=v_I.
 \eeq
 It follows from \eqref{V1} that $A_n(x)$ and  $A_p(x)$ are bounded from below, i.e.
 \beq
 A_n(x)\geq \frac{K}{2},\ \  A_p(x) \geq \frac{K}{2},\qquad\text { a.e.  } x\in\mathbb{R}^3.\non
 \eeq
Under the transformation \eqref{transform}, the right-hand side of \eqref{TAPX} are given by:
 \beq\begin{split}
 f_1(u,v,\psi_\sigma)=&\ Ku+ e^{\frac{V_n(x)}{2}}\big(\ddiv(n_\sigma\nabla\psi_\sigma)-\tilde{R}(n_\sigma, p_\sigma)\big)\\
=&\ Ku+ \nabla u\cdot \nabla\psi_\sigma-\frac{1}{2}u\nabla\psi_\sigma\cdot\nabla V_n-u^2 e^{-\frac{V_n(x)}{2}}+uv e^{-\frac{V_p(x)}{2}}+D(x)u\\
 &\quad-e^{\frac{V_n(x)}{2}}\tilde{R}(n_\sigma, p_\sigma,x),
 \end{split}\non
 \eeq
 \beq
 \begin{split} f_2(u,v,\psi_\sigma)=&\ Kv+e^{\frac{V_p(x)}{2}}\big(\ddiv(-p_\sigma\nabla\psi_\sigma)-\tilde{R}(n_\sigma, p_\sigma)\big)\\
 =&\ Kv-\nabla v\cdot \nabla\psi_\sigma+\frac{1}{2}v\nabla\psi_\sigma\cdot\nabla V_p-v^2 e^{-\frac{V_p(x)}{2}}+uv e^{-\frac{V_n(x)}{2}}-D(x)v\\
 &\quad-e^{\frac{V_p(x)}{2}}\tilde{R}(n_\sigma, p_\sigma,x).
 \end{split}\non
 \eeq

In what follows, we first prove the local well-posedness of the transformed approximate problem \eqref{TAPX}--\eqref{initial3}.
\begin{proposition}\label{localAPX}
Suppose that (H1a)--(H3) are satisfied and $u_I$, $v_I\in L^2_{+}(\mathbb{R}^3)$. Then for any $\sigma>0$, there exists $T_{\sigma}>0$ such that problem \eqref{TAPX}--\eqref{initial3} admits a unique solution $(u,v,\psi_\sigma)$ on $[0, T_{\sigma}]$, which satisfies
 \beq
 u,v\in C([0,T_{\sigma}];L^2_{+}(\mathbb{R}^3)),\qquad\nabla u,\nabla v\in L^2(0,T_{\sigma};   \mathbf{L}^2(\mathbb{R}^3)),\non
 \eeq
 \beq
 u\nabla V_n, v\nabla V_p\in L^2(0,T_{\sigma};\mathbf{L}^2(\mathbb{R}^3)).\non
 \eeq
The potential $\psi_\sigma$ is given by
 \beq
 \psi_\sigma(x,t)=\frac{1}{S_3}\int_{\mathbb{R}^3}\frac{u(y,t)e^{-\frac{V_n(y)}{2}}-v(y,t)e^{-\frac{V_p(y)}{2}}-D(y)}{|x-y|}\,dy.
 \non
 \eeq
\end{proposition}
\begin{proof}
 We consider the following auxiliary linear problem of the transformed approximated problem \eqref{TAPX}--\eqref{initial3}, such that for any
 $\tilde{u},\tilde{v}\in C([0,T];L^2(\mathbb{R}^3)),\ \nabla \tu,\nabla \tv\in L^2(0,T;\mathbf{L}^2(\mathbb{R}^3))$, $\tu\nabla V_n, \tv\nabla V_p\in L^2(0,T;\mathbf{L}^2(\mathbb{R}^3))$,
 \beq\label{ATAPX}
 (ATAP) \begin{cases}u_t-\Delta u+A_n(x)u=f_1^+(\tilde{u},\tilde{v},\tilde{\psi}_\sigma),\\
 v_t-\Delta v+A_p(x) v=f_2^+(\tilde{u},\tilde{v},\tilde{\psi}_\sigma),\\
 u(x,t)|_{t=0}=u_I,\quad v(x,t)|_{t=0}=v_I,
 \end{cases}
 \eeq
 where
 $\tilde{\psi}_\sigma$ satisfies
 \beq
 -\Delta\tilde{\psi}_\sigma= \tilde{u}e^{-\frac{V_n(x)}{2}}-\tilde{v}e^{-\frac{V_p(x)}{2}}-D(x),
 \eeq
 and the nonlinearities $f^+_1, f^+_2$ are given by
 \beq
f^+_1(\tu,\tv,\tilde{\psi}_\sigma)
=  K\tu^+ +e^{\frac{V_n(x)}{2}}\left[{\rm div}\left(\tu^+e^{-\frac{V_n(x)}{2}}\nabla \tilde{\psi}_\sigma\right)-\tilde{R}\left(\tu^+e^{-\frac{V_n(x)}{2}},\tv^+ e^{-\frac{V_p(x)}{2}},x\right)\right],\non
\eeq
\beq
 f^+_2(\tu,\tv,\tpsi_\sigma)
= K\tv^+ + e^{\frac{V_p(x)}{2}}\left[-{\rm div}\left(\tv^+ e^{-\frac{V_p(x)}{2}}\nabla
\tilde{\psi}_\sigma\right)-\tilde{R}\left(\tu^+e^{-\frac{V_n(x)}{2}},\tv^+e^{-\frac{V_p(x)}{2}},x\right)\right],\non
\eeq
\beq
 \text{with}  \quad \tu^+:=\max\{0,\tu\},\quad  \tv^+:=\max\{0,\tv\}.\non
 \eeq
 Since now the nonlinearity $\tilde{R}$ in the approximate problem \eqref{APX} satisfies the properties in Lemma \ref{lipc}, using assumptions (H1a)--(H3) we are able to prove the local well-posedness of problem \eqref{TAPX}--\eqref{initial3}
 by the contraction mapping principle as in \cite[Theorem 2.2]{WMZ}. Since the proof is the same, we only sketch it here. Denote
 \bea \Sigma_T & = & \Big\{(u,v)\in
     C\left([0,T];L^2(\mathbb{R}^3)\times L^2(\mathbb{R}^3) \right), \non\\
     && \ \ \ (\nabla u,\nabla v)\in  L^2\left(0,T; \mathbf{L}^2(\mathbb{R}^3)
     \times  \mathbf{L}^2(\mathbb{R}^3)\right)\;
     :\non\\ && \ \ \ u(0)=u_I\geq 0,\ v(0)=v_I\geq 0,
     \ \ \max_{0\leq t\leq T}\left(\|u\|_{L^2(\mathbb{R}^3)}^2
     +\|v\|_{L^2(\mathbb{R}^3)}^2\right)\leq 2M,\non\\
 &   & \ \  \ \int_0^T\left(\|\nabla u(t)\|_{\tL}^2+\|\nabla v(t)\|_{\tL}^2\right)dt
 \leq 2M, \non\\
 && \ \ \ \int_0^T\int_{\mathbb{R}^3} \left(|u(t)\nabla V_n|^2+|v(t)\nabla V_p|^2\right)dt\leq 2M. \Big\}\non
  \eea where
 \beq
 M:=\|u_I\|_{L^2(\mathbb{R}^3)}^2+\|v_I\|_{L^2(\mathbb{R}^3)}^2.\non
 \eeq
 Then we can prove that there exists a sufficiently small $T_\sigma>0$ such that the mapping
 $\mathcal{G}: (\tilde{u},\tilde{v})\mapsto (u,v)$ defined by \eqref{ATAPX} maps $\Sigma_{T_\sigma}$ to itself and is a strict contraction. Hence, the contraction principle entails that $\mathcal{G}$ has a unique fixed point in $\Sigma_{T_\sigma}$ such that $\mathcal{G}(u,v)= (u,v)$.
  Next, due to the special structure of the approximated recombination--generation rate $\tilde{R}$, using the idea in \cite{Ga}, one can show the nonnegativity of the fixed point $(u,v)$ of $\mathcal{G}$, provided that the initial data $u_I, v_I$ are nonnegative (cf. \cite[Theorem 2.2]{WMZ}). Since $ \tu^+=\tu$, $\tv^+=\tv$ if $\tu,\tv \geq 0$, we see that $f^+_i(\tu,\tv,\tilde{\psi}_\sigma)=f_i(\tu,\tv,\tilde{\psi}_\sigma)$ $(i=1,2)$ for $\tu,\tv \geq 0$.  Thus, $(u,v)$ is the local solution of problem \eqref{TAPX}--\eqref{initial3}. The details are omitted here.
\end{proof}

\begin{lemma}\label{L21}
Assume that (H1a)--(H3) are satisfied. For any $T>0$, if $n_I, p_I\in L^r(\mathbb{R}^3)\cap L^1(\mathbb{R}^3)$, $r\in\mathbb{N}$, we have
 \beq
 \|n_\sigma(t)\|_{L^s(\mathbb{R}^3)}+\|p_\sigma(t)\|_{L^s(\mathbb{R}^3)}\leq C_{T},\qquad 0\leq t\leq T,\qquad s=1,...,r.
 \eeq
Moreover,
 \beq
 \label{psibound}
 \|\nabla\psi_\sigma\|_{L^{\infty}(\mathbb{R}^3)}\leq C_T,\qquad 0\leq t\leq T,
 \eeq
provided that $n_I, p_I\in L^4(\mathbb{R}^3)\cap L^1(\mathbb{R}^3)$. In particular, the constant $C_T$ is independent of $\sigma>0$.
\end{lemma}
\begin{proof}
Integrating the equations for $n_\sigma$ and $p_\sigma$ in \eqref{APX} on $\mathbb{R}^3$, we infer from \eqref{bound}, (H2a)--(H2b) that
 \begin{eqnarray}
 & &\frac{d}{dt}\big(\|n_\sigma\|_{L^1(\mathbb{R}^3)}+\|p_\sigma\|_{L^1(\mathbb{R}^3)}\big)\non\\
&=& -2\int_{\mathbb{R}^3} R\left(\frac{n_\sigma}{1+\sigma n_\sigma},\frac{p_\sigma}{1+\sigma p_\sigma}, x\right)dx\non\\
&=& -2\int_{\mathbb{R}^3}\left(\frac{n_\sigma p_\sigma}{(1+\sigma n_\sigma)(1+\sigma p_\sigma)}-\mu_n\mu_p\right)F\left(\frac{n_\sigma}{1+\sigma n_\sigma},\frac{p_\sigma}{1+\sigma p_\sigma}\right)dx\non\\
&\leq &\ 2\int_{\mathbb{R}^3}\mu_n\mu_p F\left(\frac{n_\sigma}{1+\sigma n_\sigma},\frac{p_\sigma}{1+\sigma p_\sigma}\right)dx\non\\
&\leq &\ 2c_2\int_{\mathbb{R}^3} \mu_n\mu_p (1+n_\sigma+p_\sigma)dx\non\\
& \leq& \ C\big(1+\|n_\sigma\|_{L^1(\mathbb{R}^3)}+\|p_\sigma\|_{L^1(\mathbb{R}^3)}\big),\non
 \end{eqnarray}
which yields
 \beq
 \|n_\sigma(t)\|_{L^1(\mathbb{R}^3)}+\|p_\sigma(t)\|_{L^1(\mathbb{R}^3)}\leq (\|n_I\|_{L^1(\mathbb{R}^3)}+\|p_I\|_{L^1(\mathbb{R}^3)}+1) e^{CT},\qquad\forall \, t\in[0,T].
 \label{L1s}
 \eeq
Next, multiplying the equations for $n_\sigma$ and $p_\sigma$ by $n_\sigma^r, p_\sigma^r$ $(r\in\mathbb{N})$ respectively, integrating on $\mathbb{R}^3$ and adding the resultants, we obtain that
 \begin{eqnarray}
 && \frac{1}{r+1}\frac{d}{dt}\int_{\mathbb{R}^3}(n_\sigma^{r+1}+p_\sigma^{r+1})dx
+\frac{4r}{(r+1)^2}\int_{\mathbb{R}^3}\left(\left|\nabla\Big(n_\sigma^{\frac{r+1}{2}}\Big)\right|^2+\left|\nabla\Big(p_\sigma^{\frac{r+1}{2}}\Big)\right|^2\right)dx\non\\
&=&\frac{r}{r+1}\int_{\mathbb{R}^3}\Delta\psi_\sigma(n_\sigma^{r+1}-p_\sigma^{r+1})+\frac{r}{r+1}\int_{\mathbb{R}^3}(\Delta V_n n_\sigma^{r+1}+\Delta V_p p_\sigma^{r+1})dx
\non\\
&&-\int_{\mathbb{R}^3}R\left(\frac{n_\sigma}{1+\sigma n_\sigma},\frac{p_\sigma}{1+\sigma p_\sigma},x\right)(n_\sigma^r+p_\sigma^r)dx.\non
 \end{eqnarray}
Since $n_\sigma, p_\sigma\geq 0$, we observe that
 \beq
 \begin{split}
 &\ -\int_{\mathbb{R}^3} R\left(\frac{n_\sigma}{1+\sigma n_\sigma},\frac{p_\sigma}{1+\sigma p_\sigma},x\right)(n_\sigma^r+p_\sigma^r)dx\\
=&\ -\int_{\mathbb{R}^3}\frac{n_\sigma p_\sigma}{(1+\sigma n_\sigma)(1+\sigma p_\sigma)}F\left(\frac{n_\sigma}{1+\sigma n_\sigma},\frac{p_\sigma}{1+\sigma p_\sigma}\right)(n_\sigma^r+p_\sigma^r)dx\\
&\ +\int_{\mathbb{R}^3}\mu_n\mu_p F\left(\frac{n_\sigma}{1+\sigma n_\sigma},\frac{p_\sigma}{1+\sigma p_\sigma}\right)(n_\sigma^r+p_\sigma^r)dx\\
\leq&\ C\int_{\mathbb{R}^3}(n_\sigma^{r}+p_\sigma^{r}+n_\sigma^{r+1}+p_\sigma^{r+1})dx.
 \end{split}\non
 \eeq
On the other hand, from the Poisson equation and the elementary calculation
 \beq
  (a^{r+1}-b^{r+1})(b-a)=-\sum\limits_{k=0}^{r}a^{r-1}b^{k}(a-b)^2\leq 0,\qquad \forall a,b\geq0,\non
  \eeq
we infer from (H3) that
 \bea  \frac{r}{r+1}\int_{\mathbb{R}^3}\Delta\psi_\sigma(n_\sigma^{r+1}-p_\sigma^{r+1})dx
 &=& \frac{r}{r+1}\int_{\mathbb{R}^3}(p_\sigma-n_\sigma+D(x))(n_\sigma^{r+1}-p_\sigma^{r+1})dx\non\\
 &\leq& \|D(x)\|_{L^{\infty}(\mathbb{R}^3)}\int_{\mathbb{R}^3}(n_\sigma^{r+1}+p_\sigma^{r+1})dx.\non
 \eea
Besides, it follows from (H1a) that
 \beq
 \frac{r}{r+1}\int_{\mathbb{R}^3}(\Delta V_n n_\sigma^{r+1}+ \Delta V_p p_\sigma^{r+1})dx
 \leq K\int_{\mathbb{R}^3}(n_\sigma^{r+1}+p_\sigma^{r+1})dx.\non
 \eeq
Summing up, we have
 \bea
 &&\frac{1}{r+1}\frac{d}{dt}\int_{\mathbb{R}^3}(n_\sigma^{r+1}+p_\sigma^{r+1})dx
+\frac{4r}{(r+1)^2}\int_{\mathbb{R}^3}\left(\left|\nabla\Big(n_\sigma^{\frac{r+1}{2}}\Big)\right|^2
+\left|\nabla\Big(p_\sigma^{\frac{r+1}{2}}\Big)\right|^2\right)dx\non\\
&\leq& C\int_{\mathbb{R}^3}\left(n_\sigma^{r}+p_\sigma^{r}+n_\sigma^{r+1}+p_\sigma^{r+1}\right)dx\non\\
&\leq&  C\int_{\mathbb{R}^3}\left(1+n_\sigma^{r+1}+p_\sigma^{r+1}\right)dx.\non
 \eea
Then it follows from the Gronwall inequality that
 \beq
 \|n_\sigma\|^{r+1}_{L^{r+1}(\mathbb{R}^3)}+\|p_\sigma\|^{r+1}_{L^{r+1}(\mathbb{R}^3)}
+\int_0^T\int_{\mathbb{R}^3}\left(\left|\nabla\Big(n_\sigma^{\frac{r+1}{2}}\Big)\right|^2
+\left|\nabla\Big(p_\sigma^{\frac{r+1}{2}}\Big)\right|^2\right)dxdt\leq C_{T},\quad r\in\mathbb{N},\non
 \eeq
provided that $\int_{\mathbb{R}^3}\left(n_{I}^{r+1}+p_{I}^{r+1}\right)dx<\infty.$

Since
 \beq
 \psi_\sigma=\frac{1}{S_3}\int_{\mathbb{R}^3}\frac{n_\sigma(y,t)-p_\sigma(y,t)-D(y)}{|x-y|}dy,\non
 \eeq
 we have
 \beq
 \nabla\psi_\sigma(x,t)=\frac{1}{S_3}\int_{\mathbb{R}^3}\frac{(n_\sigma(y,t)-p_\sigma(y,t)-D(y))(x-y)}{|x-y|^3}dy.\non
 \eeq
 For any $x\in\mathbb{R}^3$,
 \beq
 |\nabla\psi_\sigma(x,t)|\leq C\int_{\mathbb{R}^3}\frac{|n_\sigma(y,t)-p_\sigma(y,t)-D(y)|}{|x-y|^2}dy.\non
 \eeq
 It follows from \cite[Corollary 2.2]{KK} (a direct consequence of the Hardy--Littlewood--Sobolev inequality \cite{Z89}) that for any $1<q'< q<\infty$ with $\frac{1}{q}=\frac{1}{q'}-\frac13$,
   \beq
   \|\nabla\psi_\sigma\|_{\mathbf{L}^q(\mathbb{R}^3)}\leq C\|n_\sigma-p_\sigma-D(x)\|_{L^{q'}(\mathbb{R}^3)}.\non
   \eeq
 Besides, if $n_\sigma,p_\sigma\in L^{\infty}(0,T;L^{q'}(\mathbb{R}^3)\cap L^1(\mathbb{R}^3))$ with $q'>3$, then we have
 \beq
 \nabla\psi_\sigma\in L^{\infty}(0,T; \mathbf{L}^{\infty}(\mathbb{R}^3)).\non
 \eeq
 The proof is complete.
\end{proof}

Based on the  \emph{a priori} estimates obtained in Lemma \ref{L21}, we can prove existence of global solutions to problem \eqref{APX}--\eqref{initial2}.

\begin{proposition}\label{globalAPX}
Suppose that all assumptions in Proposition \ref{localAPX} are satisfied.
Assume in addition that $n_I, p_I\in L^4(\mathbb{R}^3)$. Then the local solution $(n_\sigma, p_\sigma, \psi_\sigma)$ obtained in Proposition \ref{localAPX} is global.
\end{proposition}
\begin{proof}
 Multiplying the first two equations in \eqref{TAPX} by $u$ and $v$ respectively, integrating on $\mathbb{R}^3$, we get
 \bea
 && \frac{1}{2}\frac{d}{dt}\|u\|_{L^2(\mathbb{R}^3)}^2+\|\nabla u\|_{\mathbf{L}^2(\mathbb{R}^3)}^2+\int_{\mathbb{R}^3}A_n(x)u^2dx\non\\
&=&\int_{\mathbb{R}^3}f_1(u,v,\psi_\sigma)udx\non\\
&=& K\int_{\mathbb{R}^3}u^2 dx+\frac{1}{2}\int_{\mathbb{R}^3}u^2\Delta\psi_\sigma dx-\frac{1}{2}\int_{\mathbb{R}^3}u^2\nabla\psi_\sigma\cdot\nabla V_ndx\nonumber\\
&& -\int_{\mathbb{R}^3}ue^{\frac{V_n(x)}{2}}R\left(\frac{ue^{-\frac{V_n(x)}{2}}}{1+\sigma ue^{-\frac{V_n(x)}{2}}},\frac{ve^{-\frac{V_p(x)}{2}}}{1+\sigma ve^{-\frac{V_p(x)}{2}}},x\right)dx,\non
 \eea
 \bea
 &&\frac{1}{2}\frac{d}{dt}\|v\|_{L^2(\mathbb{R}^3)}^2+\|\nabla v\|_{\mathbf{L}^2(\mathbb{R}^3)}^2+\int_{\mathbb{R}^3}A_p(x)v^2dx
\non\\
&=&\int_{\mathbb{R}^3}f_2(u,v,\psi_\sigma)udx\non\\
&=&K\int_{\mathbb{R}^3}v^2 dx-\frac{1}{2}\int_{\mathbb{R}^3}v^2\Delta\psi_\sigma dx+\frac{1}{2}\int_{\mathbb{R}^3}v^2\nabla\psi_\sigma\cdot\nabla V_pdx\non\\
&& -\int_{\mathbb{R}^3}ve^{\frac{V_p(x)}{2}}R\left(\frac{ue^{-\frac{V_n(x)}{2}}}{1+\sigma ue^{-\frac{V_n(x)}{2}}},\frac{ve^{-\frac{V_p(x)}{2}}}{1+\sigma ve^{-\frac{V_p(x)}{2}}},x\right)dx.\non
 \eea
From the Poisson equation for $\psi_\sigma$, Lemma \ref{L21}, the
Gagliardo--Nirenberg inequality $$\|u\|_{L^\frac83(\mathbb{R}^3)}\leq
C\|\nabla
u\|^\frac38_{\mathbf{L}^2(\mathbb{R}^3)}\|u\|^\frac58_{L^2(\mathbb{R}^3)}$$
and the Young inequality, we have
 \bea
 && \int_{\mathbb{R}^3}(u^2-v^2)\Delta\psi_\sigma dx\non\\
 &=&\int_{\mathbb{R}^3}(u^2-v^2)(-n_\sigma+p_\sigma)dx
 +\int_{\mathbb{R}^3}D(x)(u^2-v^2)dx\non\\
 &\leq& \left(\|u\|_{L^\frac83(\mathbb{R}^3)}^2+\|p\|_{L^\frac83(\mathbb{R}^3)}^2\right)(\|n_\sigma\|_{L^4(\mathbb{R}^3)}+\|p_\sigma\|_{L^4(\mathbb{R}^3)})\non\\
 &&+ \|D(x)\|_{L^\infty(\mathbb{R}^3)}\Big(\|u\|_{L^2(\mathbb{R}^3)}^2+\|v\|_{L^2(\mathbb{R}^3)}^2\Big)\non\\
 &\leq& \frac12\left(\|\nabla u\|_{\mathbf{L}^2(\mathbb{R}^3)}^2+\|\nabla v\|_{\mathbf{L}^2(\mathbb{R}^3)}^2\right)+ C\Big(\|u\|_{L^2(\mathbb{R}^3)}^2+\|v\|_{L^2(\mathbb{R}^3)}^2\Big).\non
 \eea
We infer from the assumption $u_I, v_I\in L^2_{+}(\mathbb{R}^3)$ and the Cauchy--Schwarz inequality that
 \beq
 \int_{\mathbb{R}^3} n_I dx\leq \frac12 \int_{\mathbb{R}^3} \mu_n dx+\frac12 \int_{\mathbb{R}^3} u^2_I dx<+\infty, \quad
 \int_{\mathbb{R}^3} p_I dx\leq \frac12 \int_{\mathbb{R}^3}\mu_p dx+\frac12 \int_{\mathbb{R}^3} v^2_I dx<+\infty,\non
 \eeq
namely, $n_I,p_I\in L^1(\mathbb{R}^3)$. Lemma \ref{L21} implies that if the initial data satisfy $n_I, p_I\in L^4(\mathbb{R}^3)\cap L^1(\mathbb{R}^3)$, then the estimate \eqref{psibound} holds. As a consequence, we get
 \beq
 \begin{split}
 \left|\int_{\mathbb{R}^3}u^2\nabla\psi_\sigma\cdot\nabla V_ndx\right|
 \leq&\ \|\nabla\psi_\sigma\|_{\mathbf{L}^{\infty}(\mathbb{R}^3)}\|u\|_{L^2(\mathbb{R}^3)}\Big(\int_{\mathbb{R}^3}u^2|\nabla V_n|^2dx\Big)^{\frac{1}{2}}\\
 \leq&\ \frac{1}{8}\int_{\mathbb{R}^3}u^2|\nabla V_n|^2dx+ C_T\|u\|^2_{L^2(\mathbb{R}^3)},
 \end{split}\non
 \eeq
and
 \beq
 \begin{split}
 \left|\int_{\mathbb{R}^3}v^2\nabla\psi_\sigma\cdot\nabla V_pdx\right|
 \leq &\ \|\nabla\psi_\sigma\|_{\mathbf{L}^{\infty}(\mathbb{R}^3)}\|v\|_{L^2(\mathbb{R}^3)}\Big(\int_{\mathbb{R}^3}v^2|\nabla V_p|^2dx\Big)^{\frac{1}{2}}\\
 \leq&\ \frac{1}{8}\int_{\mathbb{R}^3}v^2|\nabla V_p|^2dx+ C_T\|v\|^2_{L^2(\mathbb{R}^3)}.
 \end{split}\non
 \eeq
It follows from the nonnegativity of $(u,v)$ and assumptions (H1b), (H2b) that
 \bea
 && -\int_{\mathbb{R}^3}\left(u e^{\frac{V_n(x)}{2}}+v e^{\frac{V_p(x)}{2}}\right)R\left(\frac{ue^{-\frac{V_n(x)}{2}}}{1+\sigma ue^{-\frac{V_n(x)}{2}}},\frac{ve^{-\frac{V_p(x)}{2}}}{1+\sigma ve^{-\frac{V_p(x)}{2}}},x\right)dx\non\\
 &\leq& \int_{\mathbb{R}^3}\left(u e^{\frac{V_n(x)}{2}}+v e^{\frac{V_p(x)}{2}}\right) F\left(\frac{ue^{-\frac{V_n(x)}{2}}}{1+\sigma ue^{-\frac{V_n(x)}{2}}},\frac{ve^{-\frac{V_p(x)}{2}}}{1+\sigma  ve^{-\frac{V_p(x)}{2}}}\right)\mu_n\mu_p dx\non\\
 &\leq& C\int_{\mathbb{R}^3} \left(u e^{\frac{V_n(x)}{2}}+v e^{\frac{V_p(x)}{2}}\right) \left(\mu_n\mu_p+u e^{-\frac{V_n(x)}{2}}+v e^{-\frac{V_p(x)}{2}} \right)dx\non\\
 &\leq& C\left(1+\|u\|^2_{L^2(\mathbb{R}^3)}+\|v\|^2_{L^2(\mathbb{R}^3)}\right).\non
 \eea
Recalling the definitions of $A_n$, $A_p$ and \eqref{V1}, we have
\bea
 \int_{\mathbb{R}^3}A_n(x)u^2dx &=&\int_{\mathbb{R}^3}\Big(\frac{1}{4}|\nabla V_n|^2-\frac{1}{2}\Delta V_n+K\Big)u^2dx\non\\
&\geq&\frac{1}{4}\int_{\mathbb{R}^3}u^2|\nabla V_n|^2dx,\eea and
\bea
 \int_{\mathbb{R}^3}A_p(x)v^2dx &=&\int_{\mathbb{R}^3}\Big(\frac{1}{4}|\nabla V_p|^2-\frac{1}{2}\Delta V_p+K\Big)v^2dx\non\\
&\geq&\frac{1}{4}\int_{\mathbb{R}^3}v^2|\nabla V_p|^2dx.\eea
As a result, we have
 \bea
 && \frac{d}{dt}\Big(\|u\|^2_{L^2(\mathbb{R}^3)}+\|v\|_{L^2(\mathbb{R}^3)}^2\Big)+\|\nabla u\|_{\mathbf{L}^2(\mathbb{R}^3)}^2+\|\nabla v\|_{\mathbf{L}^2(\mathbb{R}^3)}^2\non\\
 &&+\frac14\int_{\mathbb{R}^3}u^2|\nabla V_n|^2dx+\frac14\int_{\mathbb{R}^3}v^2|\nabla V_p|^2dx\non\\
 &\leq&
 C_T(1+\|u\|^2_{L^2(\mathbb{R}^3)}+\|v\|_{L^2(\mathbb{R}^3)}^2).\non
 \eea
 It follows from the Gronwall inequality that for any $T>0$ and $t\in[0,T]$,
 \bea
 && \|u(t)\|_{L^2(\mathbb{R}^3)}+\|v(t)\|_{L^2(\mathbb{R}^3)}\leq C_T,\\
 && \int_0^T(\|\nabla u\|_{\mathbf{L}^2(\mathbb{R}^3)}^2+\|\nabla v\|_{\mathbf{L}^2(\mathbb{R}^3)}^2)dt\leq C_T,\label{l2h1} \\
 && \int_0^T\int_{\mathbb{R}^3}(u^2|\nabla V_n|^2dx+v^2|\nabla V_p|^2)dxdt\leq C_T,\label{uvV}
 \eea
 where $C_T$ is a constant depending on $T, c_2, V_b$, $\|u_I\|_{L^2(\mathbb{R}^3)}$, $\|v_I\|_{L^2(\mathbb{R}^3)}$, $\|n_I\|_{L^4(\mathbb{R}^3)}$, $\|p_I\|_{L^4(\mathbb{R}^3)}$ but it is independent of $\sigma$.

  Then we are able to extend the local solution $(n_\sigma, p_\sigma, \psi_\sigma)$ obtained in Proposition \ref{localAPX} to the interval  $[0,T]$ for arbitrary $T>0$. The proof is complete.
\end{proof}

\textbf{Proof Theorem \ref{apx}.} Recalling the transformation \eqref{transform}, we conclude Theorem \ref{apx} from Propositions \ref{localAPX}, \ref{globalAPX}. The proof is complete.

\section{Well-posedness of the Original Problem}
\setcounter{equation}{0} In this section, we prove the existence and
uniqueness of global solutions to the Cauchy problem of original
system \eqref{DDP}. For this purpose, we shall derive some \emph{a
priori} estimates on the solutions $(n_\sigma, p_\sigma,
\psi_\sigma)$ of the approximate problem that are uniform in the
parameter $\sigma>0$. Then, we pass to the limit as $ \sigma\to 0$
to achieve our goal. In Lemma \ref{L21}, we have already shown the
uniform estimates on $\|n_\sigma(t)\|_{L^r(\mathbb{R}^3)}$,
$\|p_\sigma(t)\|_{L^r(\mathbb{R}^3)}$ on arbitrary interval $[0,T]$.
Next, we proceed to obtain uniform estimates on the $L^{\infty}$
norms of $n_\sigma$ and $p_\sigma$ via a Stampacchia-type $L^\infty$
estimation technique (cf. \cite{FI, GMT}). The following technical lemma
 plays an important role in the proof.
 It shows that a nonnegative, non-increasing function will vanish at some finite value under suitable growth condition that indicates
 certain rapid decay of the function.

\begin{lemma}\label{LM}
Suppose that $\omega(k)$ is a nonnegative non-increasing function on $[k_0, +\infty)$, and there are positive constants $\gamma,\beta$ such that
 \beq
 \omega(\hat{k})\leq M(k) (\hat{k}-k)^{-\gamma}\omega(k)^{1+\beta}, \quad \forall\, \hat{k}>k\geq k_0,\non
 \eeq
 where the function $M(k)$ is non-decreasing and satisfies
 \beq
  0\leq k^{-\gamma}M(k)\leq M_0, \quad \forall \, k\in [k_0,+\infty).\non
 \eeq
 Then
 \beq
 \omega(k^*)=0 \ \ \  \text{with}\ \ k^*=2k_0\left(1+2^{\frac{1+2\beta}{\beta^2}}M_0^{\frac{1+\beta}{\beta\gamma}}\omega(k_0)^{\frac{1+\beta}{\gamma}}\right).\non
 \eeq
\end{lemma}

\begin{remark}
Readers may refer to \cite[Lemma 2.3]{FI} for the proof of Lemma \ref{LM}. Besides, we note that the conclusion of Lemma \ref{LM} implies that  $\omega(k)=0$ for all $k\geq k^*$.
\end{remark}

\begin{lemma}\label{Linfty} Suppose that all assumptions in Theorem \ref{apx} are satisfied. Assume in addition that $n_I,p_I \in L^{\infty}(\mathbb{R}^3)$.
Then for any $T>0$, we have
 \beq
 \|n_\sigma(t)\|_{L^{\infty}(\mathbb{R}^3)}\leq C_T,\quad \|p_\sigma(t)\|_{L^{\infty}(\mathbb{R}^3)}\leq C_T,\quad 0\leq t \leq T,
 \eeq
where the constant $C_T$ is independent of $\sigma>0$.
\end{lemma}
\begin{proof}
Denote
 \beq
 k_0:=\max\{\|n_I\|_{L^{\infty}(\mathbb{R}^3)}, \|p_I\|_{L^{\infty}(\mathbb{R}^3)}\}\geq 0.\label{k0}
 \eeq
For any $k\geq k_0$, we introduce the sets
 \beq
 B_{nk}(t)= \{x\in\mathbb{R}^3, n_\sigma(x,t)>k \},\quad B_{pk}(t)=\{x\in\mathbb{R}^3, p_\sigma(x,t)>k \},\non
 \eeq
 and for arbitrary $T>0$, we set
 \beq
 \omega_T(k)=\sup\limits_{0\leq t\leq T}(|B_{nk}(t)|+|B_{pk}(t)|).\non
 \eeq
It is obvious that $\omega_T(k)$ is a nonnegative, non-increasing function on $[k_0,+\infty)$. Moreover, we infer from the $L^1$-estimate \eqref{L1s} that for arbitrary but fixed $T>0$, $\omega_T(k)$ is bounded.

For any $k\geq k_0$ and $f(\cdot)\in H^1(\mathbb{R}^3)$, it follows
from $k\geq0$ and $(f-k)^+\leq |f|$ that $(f-k)^+\in
L^2(\mathbb{R}^3)$. Moreover, for the weak derivative of
$(f-k)^+$, we have (see \cite[Lemma 7.6]{GT} or \cite[Theorem 1.56]{GMT})
 \beq \nabla (f-k)^+=\left\{
\begin{aligned}
&\nabla f \ \ \text{if }f>k\ \ &\\
&0\ \ \ \ \ \text{if }f\leq k. &\end{aligned}\right.\non
 \eeq Hence,
$(f-k)^+\in H^1(\mathbb{R}^3).$ Recalling Theorem \ref{apx} and the
lower boundedness of $V_n$ and $V_p$, we find that
 \beq n_{\sigma},
\;\ p_{\sigma}\in  L^2(0,T; H^1(\mathbb{R}^3)),
\non
 \eeq which indicates
  \beq (n_\sigma-k)^+,\;\ (p_\sigma-k)^+\in
L^2(0,T; H^1(\mathbb{R}^3)). \non
 \eeq
Now, multiplying the first and second equation in \eqref{APX} by
$(n_\sigma-k)^+$ and $(p_\sigma-k)^+$ ($k\geq k_0$), respectively, integrating on
$\mathbb{R}^3\times (0, t)$ and adding the resultants together, we obtain
 \bea
 &&  \frac{1}{2} \left(\|(n_\sigma(t)-k)^+\|^2_{L^2(\mathbb{R}^3)}+\|(p_\sigma(t)-k)^+\|^2_{L^2(\mathbb{R}^3)}\right)\non\\
 &&\ + \int_0^t \Big( \|\nabla(n_\sigma-k)^+\|^2_{\mathbf{L}^2(\mathbb{R}^3)}+\|\nabla(p_\sigma-k)^+\|^2_{\mathbf{L}^2(\mathbb{R}^3)} \Big) d\tau \non\\
 &=& \int_0^t \int_{\mathbb{R}^3} \left[\mathrm{div}(n_\sigma\nabla(\psi_\sigma+V_n))(n_\sigma-k)^+ +\mathrm{div}(p_\sigma\nabla(-\psi_\sigma+V_p))(p_\sigma-k)^+\right]dx d\tau\non\\
 & &- \int_0^t\int_{\mathbb{R}^3}R\left(\frac{n_\sigma}{1+\sigma n_\sigma},\frac{p_\sigma}{1+\sigma p_\sigma},x\right)[(n_\sigma-k)^{+}+(p_\sigma-k)^+]dxd\tau\non\\
 &=:& \int_0^t (I_1+I_2) d\tau,\label{Li1}
 \eea
 where we have used the facts that $\|(n_I-k)^+\|^2_{L^2(\mathbb{R}^3)}=\|(p_I-k)^+\|^2_{L^2(\mathbb{R}^3)}=0$ due to the assumptions $n_I, p_I \in L^\infty(\mathbb{R}^3)$ and  $k\geq k_0$ (cf. \eqref{k0}).

 Integrating by parts and using the Poisson equation for $\psi_\sigma$, we expand the term $I_1$ as follows
 \bea
 I_1&=&\int_{\mathbb{R}^3}\ddiv\left[(n_{\sigma}-k)\nabla(\psi_{\sigma}+V_n)\right](n_{\sigma}-k)^+dx+k\int_{\mathbb{R}^3}\Delta(\psi_{\sigma}+V_n)(n_{\sigma}-k)^+dx\non\\
 &&+\int_{\mathbb{R}^3}\ddiv\left[(p_{\sigma}-k)\nabla(-\psi_{\sigma}+V_p)\right](p_{\sigma}-k)^+dx+k\int_{\mathbb{R}^3}\Delta(-\psi_{\sigma}+V_p)(p_{\sigma}-k)^+dx\non\\
&=&k\int_{\mathbb{R}^3} \left[\Delta(\psi_\sigma+V_n)(n_\sigma-k)^+ +\Delta(-\psi_\sigma+V_p)(p_\sigma- k)^+\right]dx\non\\
 && -\int_{\mathbb{R}^3}(n_\sigma-k)\nabla(\psi_\sigma+V_n)\cdot\nabla (n_\sigma-k)^+ dx\non\\
 && -\int_{\mathbb{R}^3} (p_\sigma-k)\nabla(-\psi_\sigma+V_p)\cdot\nabla (p_\sigma-k)^+ dx\non\\
 &=&k\int_{\mathbb{R}^3} \Delta\psi_\sigma[(n_\sigma-k)^+-(p_\sigma- k)^+]dx
    +k\int_{\mathbb{R}^3} \left[ \Delta V_n (n_\sigma-k)^+ +\Delta V_p(p_\sigma-k)^+\right]dx\non\\
 && -\int_{\mathbb{R}^3}(n_\sigma-k)^+\nabla(\psi_\sigma+V_n)\cdot\nabla (n_\sigma-k)^+ dx\non\\
 && -\int_{\mathbb{R}^3} (p_\sigma-k)^+\nabla(-\psi_\sigma+V_p)\cdot\nabla (p_\sigma-k)^+ dx\non\\
 &=& k\int_{\mathbb{R}^3} (-n_\sigma+p_\sigma+D(x)) \left[(n_\sigma-k)^+ -(p_\sigma-k)^+\right]dx\non\\
 &&    +k\int_{\mathbb{R}^3} \left[ \Delta V_n (n_\sigma-k)^+ +\Delta V_p(p_\sigma-k)^+\right]dx \non\\
&& -\frac12\int_{\mathbb{R}^3} \nabla(\psi_\sigma+V_n)\cdot\nabla[(n_\sigma-k)^+]^2dx -\frac12\int_{\mathbb{R}^3}\nabla(-\psi_\sigma+V_p)\cdot\nabla[(p_\sigma-k)^+]^2dx
\non\\
&=& k\int_{\mathbb{R}^3} (-n_\sigma+p_\sigma) \left[(n_\sigma-k)^+ -(p_\sigma-k)^+\right]dx\non\\
&& + k\int_{\mathbb{R}^3} D(x) \left[(n_\sigma-k)^+ -(p_\sigma-k)^+\right]dx\non\\
&&   +k\int_{\mathbb{R}^3} \left( \Delta V_n (n_\sigma-k)^+ +\Delta V_p(p_\sigma-k)^+\right)dx \non\\
&& + \frac12\int_{\mathbb{R}^3}(-n_\sigma+p_\sigma) ([(n_\sigma-k)^+]^2-[(p_\sigma-k)^+]^2)dx\non\\
&& + \frac12\int_{\mathbb{R}^3} D(x)([(n_\sigma-k)^+]^2-[(p_\sigma-k)^+]^2)dx\non\\
&& +\frac12\int_{\mathbb{R}^3} (\Delta V_n[(n_\sigma-k)^+]^2+\Delta V_p[(p_\sigma-k)^+]^2)dx\non\\
&:=&J_1+...+J_6.\label{I1}
\eea
It is easy to verify that
 \beq
 J_1\leq0, \quad J_4\leq 0, \quad \forall k\geq k_0.\label{J14}
 \eeq
Besides, it follows from (H2b) that
 \beq
I_2\leq c_2\int_{\mathbb{R}^3}\mu_n\mu_p(1+n_\sigma+p_\sigma)[(n_\sigma-k)^++(p_\sigma-k)^+]dx:=J_7.\label{I2}
 \eeq
Then we infer from \eqref{Li1}--\eqref{I2} that
 \bea
 &&\frac{1}{2}\left(\|(n_\sigma(t)-k)^+\|^2_{L^2(\mathbb{R}^3)}+\|(p_\sigma(t)-k)^+\|^2_{L^2(\mathbb{R}^3)}\right)\non\\
 && \ \ +
 \int_0^t\Big(\|\nabla(n_\sigma-k)^+\|^2_{\mathbf{L}^2(\mathbb{R}^3)}+\|\nabla(p_\sigma-k)^+\|^2_{\mathbf{L}^2(\mathbb{R}^3)}\Big)d\tau\non\\
 &\leq& \int_0^t(J_2+J_3+J_5+J_6+J_7)d\tau.\non
 \eea
 By assumptions (H1b) and (H3), we have
  \bea
  && \left|\int_0^t J_5 d\tau\right|+ \left|\int_0^t J_6 d\tau\right|\non\\
  &\leq&  \frac12\|D(x)\|_{L^{\infty}(\mathbb{R}^3)}\int_0^t \|(n_\sigma-k)^+\|^2_{L^2(\mathbb{R}^3)} +\|(p_\sigma-k)^+\|^2_{L^2(\mathbb{R}^3)}d\tau\non\\
  &&+ \|\Delta V_n\|_{L^{\infty}(\mathbb{R}^3)}\int_0^t \|(n_\sigma-k)^+\|^2_{L^2(\mathbb{R}^3)} d\tau
  + \|\Delta V_p\|_{L^{\infty}(\mathbb{R}^3)}\int_0^t \|(p_\sigma-k)^+\|^2_{L^2(\mathbb{R}^3)}d\tau\non\\
  &\leq& C\int_0^t \|(n_\sigma-k)^+\|_{L^2(\mathbb{R}^3)}^2+ \|(p_\sigma-k)^+\|_{L^2(\mathbb{R}^3)}^2
 d\tau.\non
  \eea
 Let $\eta>0$ be a small constant to be chosen later. Using the H\"{o}lder inequality and Gagliardo--Nirenberg inequality, we get
 \bea
 && \left|\int_0^t J_2 d\tau\right|+ \left|\int_0^t J_3 d\tau\right|\non\\
 &\leq& k\int_0^t\int_{\mathbb{R}^3}(|D(x)|+|\Delta V_n|)(n_\sigma-k)^+   dxd\tau\non\\
 && + k\int_0^t\int_{\mathbb{R}^3}(|D(x)|+|\Delta V_p|)(p_\sigma-k)^+  dxd\tau\non\\
 &\leq& k(\|D(x)\|_{L^{\infty}(\mathbb{R}^3)}+\|\Delta V_n\|_{L^{\infty}(\mathbb{R}^3)})\int_0^t \|(n_\sigma-k)^+\|_{L^3(\mathbb{R}^3)}|B_{nk}|^{\frac{2}{3}}d\tau\non\\
 && + k(\|D(x)\|_{L^{\infty}(\mathbb{R}^3)}+\|\Delta V_p\|_{L^{\infty}(\mathbb{R}^3)})\int_0^t \|(p_\sigma-k)^+\|_{L^3(\mathbb{R}^3)}|B_{pk}|^{\frac{2}{3}}d\tau\non\\
 &\leq& C k\omega_T(k)^{\frac{2}{3}}\int_0^t\|\nabla (n_\sigma-k)^+\|_{\mathbf{L}^2(\mathbb{R}^3)}^\frac12 \|(n_\sigma-k)^+\|_{L^2(\mathbb{R}^3)}^\frac12  d\tau\non\\
 & & +C k\omega_T(k)^{\frac{2}{3}}\int_0^t\|\nabla (p_\sigma-k)^+\|_{\mathbf{L}^2(\mathbb{R}^3)}^\frac12 \|(p_\sigma-k)^+\|_{L^2(\mathbb{R}^3)}^\frac12d\tau\non\\
 &\leq& \eta T \int_0^t \|\nabla (n_\sigma-k)^+\|_{\mathbf{L}^2(\mathbb{R}^3)}^2+ \|\nabla (p_\sigma-k)^+\|_{\mathbf{L}^2(\mathbb{R}^3)}^2
 d\tau
 \non\\&& + \eta T \int_0^t \|(n_\sigma-k)^+\|_{L^2(\mathbb{R}^3)}^2+ \|(p_\sigma-k)^+\|_{L^2(\mathbb{R}^3)}^2
 d\tau\non\\
 && +\frac{C}{\eta}k^2\omega_T(k)^{\frac{4}{3}}\non
  \eea
 and
 \bea
 &&\left|\int_0^t J_7 d\tau\right|\non\\
 &\leq&\int_0^t\int_{\mathbb{R}^3}\mu_n\mu_p [(n_\sigma-k)^++(p_\sigma-k)^+ +2k+1]\left[(n_\sigma-k)^++(p_\sigma-k)^+\right]dxd\tau
 \non\\
&\leq&C\int_0^t\|(n_\sigma-k)^+\|_{L^2(\mathbb{R}^3)}^2+ \|(p_\sigma-k)^+\|_{L^2(\mathbb{R}^3)}^2d\tau\non\\
&& +C(2k+1)\int_0^t\int_{\mathbb{R}^3}(n_\sigma-k)^+dxd\tau+C(2k+1)\int_0^t\int_{\mathbb{R}^3}(p_\sigma-k)^+dxd\tau\non\\
&\leq& C\int_0^t(\|(n_\sigma-k)^+\|^2_{L^2(\mathbb{R}^3)}+\|(p_\sigma-k)^+\|^2_{L^2(\mathbb{R}^3)})d\tau\non\\
&& +C(2k+1)\omega_T(k)^{\frac{2}{3}}\int_0^t\left(\|(n_\sigma-k)^+\|_{L^3(\mathbb{R}^3)}+\|(p_\sigma-k)^+\|_{L^3(\mathbb{R}^3)}\right)d\tau\non\\
&\leq& \eta T \int_0^t \|\nabla (n_\sigma-k)^+\|_{\mathbf{L}^2(\mathbb{R}^3)}^2+\|\nabla (p_\sigma-k)^+\|_{\mathbf{L}^2(\mathbb{R}^3)}^2
 d\tau\non\\
 && + (C+\eta T) \int_0^t \|(n_\sigma-k)^+\|_{L^2(\mathbb{R}^3)}^2 +\|(p_\sigma-k)^+\|_{L^2(\mathbb{R}^3)}^2d\tau\non\\
 && +\frac{C}{\eta}(2k+1)^2\omega_T(k)^{\frac{4}{3}}.\non
\eea
Taking $$\eta=\frac{1}{2T}$$ in the above estimates, we obtain that
 \beq
 \begin{split}&
 \|(n_\sigma(t)-k)^+\|^2_{L^2(\mathbb{R}^3)}+\|(p_\sigma(t)-k)^+\|^2_{L^2(\mathbb{R}^3)}\\
 &\ +
\int_0^t(\|\nabla(n_\sigma-k)^+\|^2_{\mathbf{L}^2(\mathbb{R}^3)}+\|\nabla(p_\sigma-k)^+\|^2_{\mathbf{L}^2(\mathbb{R}^3)})d\tau\\
\leq&\ C_1\int_0^t\left(\|(n_\sigma(\tau)-k)^+\|^2_{L^2(\mathbb{R}^3)}+\|(p_\sigma(\tau)-k)^+\|^2_{L^2(\mathbb{R}^3)}\right)d\tau\\
&\ +TC_2(k^2+1)\omega_T(k)^{\frac{4}{3}}.
 \end{split}\non
 \eeq
It follows from the Gronwall inequality that for $t\in [0,T]$
 \beq
 \begin{split}
 &\|(n_\sigma(t)-k)^+\|^2_{L^2(\mathbb{R}^3)}+\|(p_\sigma(t)-k)^+\|^2_{L^2(\mathbb{R}^3)}
 \leq e^{C_1T}TC_2(k^2+1)\omega_T(k)^{\frac{4}{3}}.
 \end{split}\label{W1}
 \eeq
On the other hand, for any $t\in[0,T]$ and $\hat{k}>k\geq k_0$,
 \bea
 \|(n_\sigma(t)-k)^+\|^2_{L^2(\mathbb{R}^3)}&=&\int_{B_{nk}(t)}|(n_\sigma(t,\cdot)-k)^+|^2dx\non\\
 & \geq & \int_{B_{n\hat{k}}(t)}(n_\sigma(t,\cdot)-k)^2dx\non\\
 &\geq&(\hat{k}-k)^2|B_{n\hat{k}}(t)|,\label{W2}
  \eea
  \bea
  \|(p_\sigma(t)-k)^+\|^2_{L^2(\mathbb{R}^3)}&=&\int_{B_{pk}(t)}|(p_\sigma(t,\cdot)-k)^+|^2dx\non\\
  &\geq & \int_{B_{p\hat{k}}(t)}(p_\sigma(t,\cdot)-k)^2dx\non\\
  &\geq&(\hat{k}-k)^2|B_{p\hat{k}}(t)|.\label{W3}
  \eea
We deduce from \eqref{W1}--\eqref{W3} that
 \beq
 \omega_T(\hat{k})\leq e^{C_1T}TC_2(k^2+1)(\hat{k}-k)^{-2}\omega_T(k)^{\frac{4}{3}},
 \quad\forall \hat{k}>k\geq k_0.\non
 \eeq
 Now in Lemma \ref{LM}, we set
 \beq
 M(k)=e^{C_1T}TC_2(k^2+1)\geq 0, \quad \gamma=2,\quad  \beta=\frac13.\non
 \eeq
 The function $M(k)$ has the following property
 \beq
 \frac{M(k)}{k^2}=\frac{k^2+1}{k^2} e^{C_1T}TC_2\leq \left(1+\frac{1}{k_0^2}\right)e^{C_1T}TC_2:=M_0, \quad \forall\ k\in [k_0,+\infty).\non
 \eeq
 Therefore, there exists a constant
 \beq
  k^*=2k_0\left(1+2^{15}M_0^2\omega_T(k_0)^\frac{2}{3}\right)>k_0,\non
 \eeq
 which is independent of $\sigma$ such that
 \beq
 \omega_T(k^*)=0.\non
 \eeq
Namely,
 \beq
 n_\sigma(x,t)\leq k^* \quad \text{   and   } \quad p_\sigma(x,t)\leq k^*,\quad
 \text{ a.e.  in }\  \mathbb{R}^3\times[0,T].\non
 \eeq
 The proof is complete.
\end{proof}

{\bf Proof of Theorem \ref{Global}.}
It follows from Lemma \ref{L21}, \eqref{l2h1} and Lemma \ref{Linfty} that the following uniform estimates (independent of the parameter $\sigma>0$) hold:
 \beq
 \|n_\sigma\|_{L^\infty(0,T; L^q(\mathbb{R}^3))}+\|p_\sigma\|_{L^\infty(0,T;L^q(\mathbb{R}^3))}\leq C_T,\label{b1}
 \eeq
 \beq
 \|n_\sigma\|_{L^2(0,T;H^1(\mathbb{R}^3))}+\|p_\sigma\|_{L^2(0,T;H^1(\mathbb{R}^3))} \leq C_T,\label{b2}
 \eeq
 \beq
 \|\nabla\psi_\sigma\|_{L^{\infty}(0,T;\mathbf{L}^{q'}(\mathbb{R}^3))}+\|\Delta\psi_\sigma\|_{L^{\infty}(0,T;L^q(\mathbb{R}^3))}\leq C_T,\label{b3}
 \eeq
 where $q\in [1,+\infty]$, $q'\in [\frac{3}{2},+\infty]$.

 Besides, we infer from \eqref{uvV}, \eqref{transform} and (H1a) that
 \bea
  && \int_0^T\int_{\mathbb{R}^3}(n_\sigma^2|\nabla V_n|^2dx+p_\sigma^2|\nabla V_p|^2)dxdt\non\\
   &\leq & e^{-V_b}\int_0^T\int_{\mathbb{R}^3}(u^2|\nabla V_n|^2dx+v^2|\nabla V_p|^2)dxdt\non\\
   &\leq& C_T.\label{b4}
 \eea
 Then by the equations for $n_\sigma$ and $p_\sigma$ in \eqref{APX} and \eqref{b1}--\eqref{b4}, we obtain that
 \beq
 \|\partial_t n_\sigma\|_{L^2(0,T; (H^1(\mathbb{R}^3))')}+\|\partial_t p_\sigma\|_{L^2(0,T; (H^1(\mathbb{R}^3))')}\leq C_T,\non
 \eeq
From the uniform estimates \eqref{b1}--\eqref{b3}, we deduce that there exist
 \beq
 n, p \in L^{\infty}(0,T;L^{\infty}(\mathbb{R}^3)),\non
 \eeq
 with
 \beq
 \nabla u, \nabla p \in L^2(0,T;\mathbf{L}^2(\mathbb{R}^3)),\quad \partial_t u, \partial_t p \in L^2(0,T; (H^1(\mathbb{R}^3))'),\non
 \eeq
 and $\psi$ with $\Delta \psi\in L^\infty(0,T; L^\infty(\mathbb{R}^3)), \nabla \psi \in L^2(0,T; \mathbf{L}^2(\mathbb{R}^3))$ such that for a sequence $\{\sigma_j\}\searrow 0$ as $j\to +\infty$ (not relabeled when taking a subsequence),
 \bea
 && n_{\sigma_j}\rightharpoonup n,\quad p_{\sigma_j}\rightharpoonup p, \ \ \text{ weakly-star in  } L^{\infty}(0,T;L^{\infty}(\mathbb{R}^3)),\non\\
 && \nabla n_{\sigma_j}\rightharpoonup\nabla n,\quad \nabla p_{\sigma_j}\rightharpoonup\nabla p,\ \ \text{ weakly in } L^2(0,T;\mathbf{L}^2(\mathbb{R}^3)),\non\\
 && \partial_t n_{\sigma_j}\rightharpoonup\partial_t n,\quad \partial_t p_{\sigma_j}\rightharpoonup\partial_t p,\ \ \text{ weakly in } L^2(0,T; (H^1(\mathbb{R}^3))'),\non\\
 && \Delta \psi_{\sigma_j}\rightharpoonup\Delta\psi,\ \  \text{ weakly-star in } L^\infty(0,T;L^\infty(\mathbb{R}^3)),\non\\
 && \nabla \psi_{\sigma_j}\rightharpoonup \nabla \psi, \ \ \text{ weakly in } L^2(0,T; \mathbf{L}^2(\mathbb{R}^3)).\non
 \eea
 Moreover, on account of the compact embedding theorem we have (up to a subsequence if necessary and without relabelling for the sake of simplicity)
 \beq
 n_{\sigma_j}\rightarrow n,\quad p_{\sigma_j}\rightarrow p,\ \ \text{strongly  in } L^2(0,T;L^2_{loc}(\mathbb{R}^3)), \quad \text{thus also a.e. in }  \mathbb{R}^3 \times (0, T).\non
 \eeq
 Then, for any $\varphi\in L^2(0,T; C_c^{\infty}(\mathbb{R}^3))$, we have
 \beq
 \begin{split}
 & \left|\int_0^T\int_{\mathbb{R}^3}(n_{\sigma_j}\nabla\psi_{\sigma_j}-n\nabla\psi)\nabla\varphi dxdt\right|\\
 \leq& \ \left|\int_0^T\int_{\mathbb{R}^3}(n_{\sigma_j}-n)\nabla\psi_{\sigma_j}\nabla\varphi dxdt\right| + \left|\int_0^T\int_{\mathbb{R}^3}(\nabla\psi_{\sigma_j}-\nabla\psi)n\nabla\varphi dx dt\right|\\
 \leq &\ \|\nabla \psi_{\sigma_j}\|_{L^\infty(0,T; L^2(\mathbb{R}^3))}\|n_{\sigma_j}- n\|_{L^2(0,T;L^2(\mathrm{supp}\varphi))}\|\nabla\varphi\|_{L^2(0,T;L^\infty(\mathbb{R}^3))} \\
 &\ + \left|\int_0^T\int_{\mathbb{R}^3}(\nabla\psi_{\sigma_j}-\nabla\psi)(n\nabla\varphi) dx dt\right|\\
 \rightarrow &\ 0, \quad \text{as} \ \sigma_j\to 0,
  \end{split}\non
 \eeq
 similarly,
  \beq
 \left|\int_0^T\int_{\mathbb{R}^3}(p_{\sigma_j}\nabla\psi_{\sigma_j}-p\nabla\psi)\nabla\varphi dxdt\right|
 \rightarrow \ 0, \quad \text{as} \ \sigma_j\to 0.\non
 \eeq
 Next, we study the convergence of the recombination-generation rate. The uniform bound \eqref{b1} and (H2b) yield that
 \beq
 \left\|\tilde{R}(n_\sigma, p_\sigma,x)\right\|_{L^2(0,T; L^2(\mathbb{R}^3))}\leq C_T.\non
 \eeq
 Thus there exists $G\in L^2(0,T; L^2(\mathbb{R}^3))$ such that (up to a subsequence)
 \beq
 \tilde{R}(n_\sigma, p_\sigma,x)\rightarrow G,\ \ \text{weakly  in } L^2(0,T; L^2(\mathbb{R}^3))\quad \text{as} \ \sigma_j\to 0.\label{cR1}
 \eeq
 For any bounded domain $\Omega\subset \mathbb{R}^3$, there holds
 \bea
 \int_0^T\int_\Omega \left(\frac{n_{\sigma_j}}{1+\sigma_jn_{\sigma_j}}-n\right)^2dxdt
 &\leq& \int_0^T\int_\Omega(n_{\sigma_j}-n)^2 dxdt +\int_0^T\int_\Omega (\sigma_jn_{\sigma_j}n)^2dxdt\non\\
 &\to& 0, \quad \text{as}\ \sigma_j\to 0,\non
 \eea
which implies that
 \beq
 \frac{n_{\sigma_j}}{1+\sigma_j n_{\sigma_j}}\rightarrow n \quad  \text{ strongly in } L^2(0,T;L^2_{loc}(\mathbb{R}^3)).\label{cn}
 \eeq
 In the same manner, we have
 \beq
 \frac{p_{\sigma_j}}{1+\sigma_j p_{\sigma_j}}\rightarrow p \quad  \text{ strongly in } L^2(0,T;L^2_{loc}(\mathbb{R}^3)).\label{cp}
 \eeq
 Since $F$ is Lip-continuous (see (H2b)), we infer from \eqref{cn} and \eqref{cp} that on any bounded domain $\Omega\subset \mathbb{R}^3$, the following convergence (up to a subsequence)
 \bea
 && \int_0^T \int_{\Omega} \left[F\left(\frac{n_{\sigma_j}}{1+\sigma_j n_{\sigma_j}},\frac{p_{\sigma_j}}{1+\sigma_j p_{\sigma_j}}\right)-F(n,p)\right]^2dxdt\non\\
 &\leq& C \int_0^T \int_{\Omega} \left(\frac{n_{\sigma_j}}{1+\sigma_j n_{\sigma_j}}-n\right)^2+\left(\frac{p_{\sigma_j}}{1+\sigma_j p_{\sigma_j}}-p\right)^2 dxdt\non\\
 &\rightarrow&  0, \quad \text{as} \ \sigma_j\to 0,\non
 \eea
 namely,
\bea
 && F\left(\frac{n_{\sigma_j}}{1+\sigma_j n_{\sigma_j}},\frac{p_{\sigma_j}}{1+\sigma_j p_{\sigma_j}}\right)\rightarrow F(n,p),\non\\
  && \qquad\ \  \text{ strongly in } L^2(0,T;L^2_{loc}(\mathbb{R}^3))\  \text{and a.e. in } \mathbb{R}^3 \times (0, T).\label{cF}
 \eea
 As a result, we have the point-wise convergence of $\tilde{R}$ (up to a subsequence)
 \beq
  \tilde{R}(n_\sigma, p_\sigma,x)\rightarrow R(n,p,x),\quad \text{a.e. in }  \mathbb{R}^3 \times (0, T),\label{cR2}
 \eeq
 which together with \eqref{cR1} implies that $G=R(n,p,x)$ and
 \beq
 \tilde{R}(n_\sigma, p_\sigma,x)\rightharpoonup R(n,p,x),\quad \text{weakly  in } L^2(0,T; L^2(\mathbb{R}^3))\quad \text{as} \ \sigma_j\to 0.\non
 \eeq

 Based on the above convergent results, now we are able to pass to the limit by letting $\sigma_j\to 0$ in the approximate problem \eqref{APX}
and obtain a global weak solution $(n,p,\psi)$ of problem \eqref{DDP}--\eqref{initial}. The system \eqref{DDP} is satisfied in the following sense that for any
$\varphi\in C^\infty_0([0,T)\times \mathbb{R}^3)$,
 \bea
 && -\int_0^T\int_{\mathbb{R}^3} n \varphi_t dx dt+\int_0^T\int_{\mathbb{R}^3}(\nabla n+n\nabla(\psi+V_n))\cdot \nabla\varphi dxdt \non\\
 && \qquad +\int_0^T\int_{\mathbb{R}^3}R(n,p,x)\varphi dxdt=\int_{\mathbb{R}^3} n_I \varphi (x,0) dx, \label{ex1}
 \eea
 \bea
  && -\int_0^T\int_{\mathbb{R}^3} p \varphi_t dx dt+\int_0^T\int_{\mathbb{R}^3}(\nabla p+p\nabla(-\psi+V_p))\cdot \nabla\varphi dxdt\non\\
  && \qquad +\int_0^T\int_{\mathbb{R}^3}R(n,p,x)\varphi dxdt=\int_{\mathbb{R}^3} p_I \varphi (x, 0) dx,\label{ex2}
  \eea
  \beq
  \label{ex3}
  \int_0^T\int_{\mathbb{R}^3}\nabla\psi\cdot\nabla\varphi dxdt=\int_0^T\int_{\mathbb{R}^3}(n-p-D)\varphi dxdt.
  \eeq

Finally, we prove the uniqueness of global solutions to problem \eqref{DDP}--\eqref{initial}.
Let $(n_i,p_i,\psi_i)$ $(i=1,2)$ be two solutions to problem \eqref{DDP}--\eqref{initial} with initial data $n_{Ii}, p_{Ii}$. Set now
$$n=n_1-n_2, \quad  p=p_1-p_2, \quad \psi=\psi_1-\psi_2,\quad n_I=n_{I1}-n_{I2}, \quad  p_I=p_{I1}-p_{I2}.$$
Taking the difference of the equations for $n_1$ and $n_2$, testing the resultant by $n$, we find that
 \bea
 &&\frac12\|n(t)\|^2_{L^2(\mathbb{R}^3)}+\int_0^t\|\nabla n\|^2_{\mathbf{L}^2(\mathbb{R}^3)}d\tau\non\\
 &=&\frac12\|n_I\|^2_{L^2(\mathbb{R}^3)}-\int_0^t\int_{\mathbb{R}^3}n\nabla n\cdot\nabla V_ndxd\tau
 -\int_0^t\int_{\mathbb{R}^3}\left(n\nabla\psi_1+n_2\nabla\psi\right)\nabla ndxd\tau\non\\
 &&-\int_0^t\int_{\mathbb{R}^3}\Big(R(n_1, p_1,x)-R(n_2,p_2,x)\Big)ndxd\tau\non\\
 &:=& \frac12\|n_I\|^2_{L^2(\mathbb{R}^3)}+E_1+E_2+E_3.\label{d1}
 \eea
 Using the uniform estimates for $n_i, p_i, \psi_i$, we have
 \bea
 |E_1+E_2| &\leq&\frac{1}{2}\Big(\|\Delta V_n\|_{L^{\infty}(\mathbb{R}^3)}+\|\Delta\psi_1\|_{L^{\infty}(\mathbb{R}^3)}\Big)\int_0^t\|n\|^2_{L^2(\mathbb{R}^3)}d\tau\non\\
 &&+ \|n_2\|_{L^\infty(0,T;L^{3}(\mathbb{R}^3))}\int_0^t\left(\e\|\nabla n\|^2_{\mathbf{L}^2(\mathbb{R}^3)}+C_{\e}\|\nabla\psi\|^2_{\mathbf{L}^6(\mathbb{R}^3)}\right)d\tau.\label{d3}
  \eea
  On the other hand,  by (H2a)--(H2b), we deduce that
  \bea
  |E_3|&\leq & \left|\int_0^t\int_{\mathbb{R}^3} (F(n_1,p_1)-F(n_2,p_2))n\mu_n\mu_p dxd\tau\right|\non\\
  &&+ \left|\int_0^t\int_{\mathbb{R}^3} (F(n_1,p_1)-F(n_2,p_2))n_1p_1 ndxd\tau\right|\non\\
  &&+ \left|\int_0^t\int_{\mathbb{R}^3} F(n_2,p_2)(n_1 p+np_2)n dxd\tau\right|\non\\
  &\leq& C\left(V_b,\|n_i\|_{L^\infty(0,T;L^{\infty}(\mathbb{R}^3))},\|p_i\|_{L^\infty(0,T;L^{\infty}(\mathbb{R}^3))}\right)\non\\
  &&\times \int_0^t(\|n\|^2_{L^2(\mathbb{R}^3)}+\|p\|^2_{L^2(\mathbb{R}^3)})d\tau.\label{d4}
 \eea
 In a similar manner, we have the following estimate for $p$
  \bea
 && \frac12\|p\|^2_{L^2(\mathbb{R}^3)}+\int_0^t\|\nabla p\|^2_{\mathbf{L}^2(\mathbb{R}^3)}d\tau\non\\
 &\leq&\frac12\|p_I\|^2_{L^2(\mathbb{R}^3)}+\frac{1}{2}\Big(\|\Delta V_p\|_{L^{\infty}(\mathbb{R}^3)}+\|\Delta \psi_1\|_{L^{\infty}(\mathbb{R}^3)}\Big)\int_0^t\|p\|^2_{L^2(\mathbb{R}^3)}d\tau\non\\
 &&+\|p_2\|_{L^\infty(0,T;L^{3}(\mathbb{R}^3))}\int_0^t\left(\e\|\nabla p\|^2_{\mathbf{L}^2(\mathbb{R}^3)}+C_{\e}\|\nabla\psi\|^2_{\mathbf{L}^6(\mathbb{R}^3)}\right)d\tau\non\\
 &&+C\left(V_b,\|n_i\|_{L^\infty(0,T;L^{\infty}(\mathbb{R}^3))},\|p_i\|_{L^\infty(0,T;L^{\infty}(\mathbb{R}^3))}\right)  \non\\
 && \quad \times \int_0^t(\|n\|^2_{L^2(\mathbb{R}^3)}+\|p\|^2_{L^2(\mathbb{R}^3)})d\tau.
 \label{d5}
 \eea
 Since $\psi$ satisfies the Poisson equation $-\Delta \psi=n-p$, then it follows from \cite[Corollary 2.2]{KK} that
 \beq
 \|\nabla\psi\|^2_{\mathbf{L}^6(\mathbb{R}^3)}\leq C\left(\|n\|^2_{L^2(\mathbb{R}^3)}+\|p\|^2_{L^2(\mathbb{R}^3)}\right). \label{d6}
 \eeq
Therefore, taking $\e$ sufficiently small satisfying
 \beq
 0<\e\leq \frac12\min\left\{1, \|n_2\|_{L^\infty(0,T;L^3(\mathbb{R}^3))}^{-1}, \|p_2\|_{L^\infty(0,T;L^3(\mathbb{R}^3))}^{-1}  \right\},\non
 \eeq
 we deduce from \eqref{d1}--\eqref{d6} that
 \beq
 \begin{split}&
 \|n(t)\|^2_{L^2(\mathbb{R}^3)}+\|p(t)\|^2_{L^2(\mathbb{R}^3)}+\int_0^t\left(\|\nabla n\|^2_{\mathbf{L}^2(\mathbb{R}^3)}+\|\nabla p\|^2_{\mathbf{L}^2(\mathbb{R}^3)}\right)d\tau\\
 \leq&\  \|n_I\|^2_{L^2(\mathbb{R}^3)}+\|p_I\|^2_{L^2(\mathbb{R}^3)}+
 C_T\int_0^t\left(\|n\|^2_{L^2(\mathbb{R}^3)}+\|p\|^2_{L^2(\mathbb{R}^3)}\right)d\tau.
  \end{split}\non
  \eeq
From the Gronwall inequality, we can conclude the continuous dependence on the initial data that
 \beq
 \|n(t)\|^2_{L^2(\mathbb{R}^3)}+\|p(t)\|^2_{L^2(\mathbb{R}^3)}\leq \left(\|n_I\|^2_{L^2(\mathbb{R}^3)}+\|p_I\|^2_{L^2(\mathbb{R}^3)}\right) e^{C_T t},\quad\forall \, t\in[0,T],\non
 \eeq
 which yields the uniqueness. The proof is complete.

\section{Long-time Behavior}\setcounter{equation}{0}

In section 3 we have proved the existence and uniqueness of global solutions to problem \eqref{DDP}--\eqref{initial}. However, the global-in-time estimates for the solution $(n,p,\psi)$ depends on $T$ that can be chosen arbitrary. In this section, we extend the results in \cite{WMZ,FFM} to our current problem \eqref{DDP}--\eqref{initial}. For this purpose, we first need to obtain some uniform-in-time estimates on the global solution.

Let $\alpha=\int_{\mathbb{R}^3} (n_I-p_I) dx$. We easily see from \eqref{DDP} that the difference of charges is conserved for all $t>0$:
 \beq
  \int_{\mathbb{R}^3} (n(t,\cdot) -p(t,\cdot)) dx=\alpha.\label{con}
 \eeq
The relative entropy associated with \eqref{DDP} is as follows:
 \beq
 \begin{split}
 e(t):=&\int_{\mathbb{R}^3}\left[n\left(\ln\frac{n}{n_{\infty}}-1\right)+n_{\infty}\right]dx
 +\int_{\mathbb{R}^3}\left[p\left(\ln\frac{p}{p_{\infty}}-1\right)+p_{\infty}\right]dx\\
 &+\frac{1}{2}\int_{\mathbb{R}^3}|\nabla\psi-\nabla\psi_{\infty}|^2dx,
 \end{split}\non
 \eeq
 where $(n_{\infty},p_{\infty},\psi_{\infty})$ is the steady state of system \eqref{DDP} that
 satisfies
 \beq
 \begin{cases}
 n_{\infty}(x)=D_n e^{-\psi_{\infty}}\mu_n,\quad D_n\in\mathbb{R}^+,\\
 p_{\infty}(x)=D_p e^{\psi_{\infty}}\mu_p, \quad D_n\in\mathbb{R}^+,\\
 n_{\infty}p_{\infty}=\mu_n\mu_p,\quad \int_{\mathbb{R}^3}n_{\infty}dx-\int_{\mathbb{R}^3}p_{\infty}dx=\alpha,\\
 -\Delta\psi_{\infty}=n_{\infty}-p_{\infty}-D(x).
 \end{cases}\label{sta}
 \eeq
 \begin{remark}
 Denote $I=\int_{\mathbb{R}^3} e^{-\psi_\infty-V_n(x)}dx, J=\int_{\mathbb{R}^3}
  e^{\psi_\infty-V_p(x)}dx$. Then the coefficients $D_n$ and $D_p$ in \eqref{sta} are given by (cf. \cite[Lemma 3.1]{WMZ})
 \beq
 D_n=\frac{\alpha+\sqrt{\alpha^2+4IJ}}{2I},\quad D_p=\frac{\sqrt{\alpha^2+4IJ}-\alpha}{2J}\nonumber
 \eeq
 satisfying $D_nD_p=1$.
 \end{remark}
Following the argument in \cite[Theorem 3.1]{WMZ}, where the special case $V_n=V_p=V$ was considered, we can still prove the existence and uniqueness of $(n_\infty,p_\infty,\psi_\infty)$.
 \begin{proposition}
 Suppose that assumptions (H1a), (H1b) and (H3) are satisfied. Then the stationary problem \eqref{sta} admits a unique solution $(\psi_\infty, n_\infty, p_\infty)$ such that
 $$\psi_\infty\in D^{1,2}(\mathbb{R}^3)=\left\{\phi(x) \in
 L^6(\mathbb{R}^3)\left|\right. \nabla \phi\in
\mathbf{L}^2(\mathbb{R}^3)\right\}.$$ Moreover, $\psi_\infty\in L^\infty(\mathbb{R}^3)$ and
 $\nabla \psi_\infty\in \mathbf{L}^\infty(\mathbb{R}^3)\cap \mathbf{L}^2(\mathbb{R}^3)$.
 \end{proposition}
 Then we have
\begin{lemma} \label{unieee}
Suppose the assumptions of Theorem \ref{Global} are satisfied. The global solution $(n,p)$ of problem \eqref{DDP}--\eqref{initial} satisfies
 \beq
 \sup\limits_{t\geq0}\left[\|n(t)\|_{L^r(\mathbb{R}^3)}+\|p(t)\|_{L^r(\mathbb{R}^3)}\right]<\infty,\quad\forall\, r\in[1,+\infty].\non
 \eeq
 \end{lemma}

\begin{proof}

   By a straightforward calculation, we have the dissipation of the relative entropy
   \beq
   \label{dissipation}
   \frac{d}{dt}e(t)=-\mathcal {D}(t)\leq 0,
    \eeq
    with the entropy dissipation
    \begin{eqnarray}
    \mathcal{D}(t) &=& -\int_{\mathbb{R}^3}n\left|\nabla\ln\left(\frac{n}{N}\right)\right|^2dx
     -\int_{\mathbb{R}^3}p\left|\nabla\ln\left(\frac{p}{P}\right)\right|^2dx\non\\
     &&\quad -\int_{\mathbb{R}^3}R(n,p,x)\ln\left(\frac{np}{\mu_n\mu_p}\right)dx,\non
     \end{eqnarray}
   \beq
   \text{where}\ \ N=D_n e^{-\psi(t)}\mu_n,\quad P=D_pe^{\psi(t)}\mu_p.\non
    \eeq
Based on the entropy dissipation inequality, we can obtain uniform $L^r$ bounds ($r\in [1,+\infty)$) for $n(t)$ and $p(t)$ exactly as in \cite[Lemma 4.1]{WMZ}. It only remains to show the uniform $L^\infty$ estimate. We note that $L^\infty$ bounds of solutions to a simplified drift-diffusion system (without self-consistent potential $\psi$ and with a recombination-generation rate of Shockley--Read--Hall type) have been obtained in \cite{FFM} via a Nash--Moser type iteration method and the results could be extended to the case with self-consistent potential \cite{DFFM}. For the convenience of the readers, we sketch the proof for our present case with a more general recombination-generation rate.

For $r\geq 2$, using  integration by parts and the nonnegativity of $n, p$, we get
 \bea
 &&\frac{d}{dt}\int_{\mathbb{R}^3}(n^{r+1}+p^{r+1})dx
 +\frac{4r}{r+1}\int_{\mathbb{R}^3}\left(\left|\nabla\Big(n^{\frac{r+1}{2}}\Big)\right|^2+\left|\nabla\Big(p^{\frac{r+1}{2}}\Big)\right|^2\right)dx\non\\
 &=&r\int_{\mathbb{R}^3}\Delta\psi(n^{r+1}-p^{r+1})+r\int_{\mathbb{R}^3}(\Delta V_n n^{r+1}+\Delta V_p p^{r+1})dx\non\\
 && -(r+1)\int_{\mathbb{R}^3}R(n,p,x)(n^r+p^r)dx\non\\
 &=&-r\int_{\mathbb{R}^3}(n-p)(n^{r+1}-p^{r+1})dx+r\int_{\mathbb{R}^3}D(x)(n^{r+1}-p^{r+1})dx\non\\
 && +r\int_{\mathbb{R}^3}(\Delta V_n n^{r+1}+\Delta V_p p^{r+1})dx\non\\
 &&-(r+1)\int_{\mathbb{R}^3} F(n,p)(n^{r+1}p+np^{r+1})dx+(r+1)\int_{\mathbb{R}^3}\mu^2 F(n,p)(n^r+p^r)dx\non\\
 &\leq&r(\|D\|_{L^{\infty}(\mathbb{R}^3)}+\|\Delta V_n\|_{L^{\infty}(\mathbb{R}^3)})\int_{\mathbb{R}^3}n^{r+1}dx\non\\
 && +r(\|D\|_{L^{\infty}(\mathbb{R}^3)}+\|\Delta V_p\|_{L^{\infty}(\mathbb{R}^3)})\int_{\mathbb{R}^3}p^{r+1}dx\non\\
 && +C(r+1)\int_{\mathbb{R}^3}(1+n+p)(n^r+p^r)dx\non\\
 &\leq & C(r+1)\int_{\mathbb{R}^3}(n^{r+1}+p^{r+1})dx+C(r+1)\int_{\mathbb{R}^3}(n^r+p^r)dx,\label{r1}
 \eea
 where we use the facts that
 \beq
 (n-p)(n^{r+1}-p^{r+1})\geq 0, \quad F(n,p)(n^{r+1}p+np^{r+1})\geq 0.\non
 \eeq
 Since
 \beq
 (r-1) r^\frac{2-r}{r-1}> 0\ \text{for} \ r\geq 2, \quad \lim_{r\to+\infty}(r-1) r^\frac{2-r}{r-1}=1,\non
 \eeq
 by the Young's inequality and the uniform $L^1$ estimate of $n,p$, we deduce that
 \bea
 \int_{\mathbb{R}^3}(n^r+p^r)dx
 &\leq&
 \frac{1}{r^2}\int_{\mathbb{R}^3}(n+p)dx+(r-1) r^\frac{2-r}{r-1}\int_{\mathbb{R}^3}(n^{r+1}+p^{r+1})dx\non\\
 &\leq&\frac{C}{r^2}+C'\int_{\mathbb{R}^3}(n^{r+1}+p^{r+1})dx,\quad \forall\, r\geq 2,\label{r2}
 \eea
 where the constants $C, C'$ are independent of $r$.
Therefore, it follows from \eqref{r1} and \eqref{r2} that
 \beq
 \begin{split}
 &\frac{d}{dt}\int_{\mathbb{R}^3}(n^{r+1}+p^{r+1})dx
 +\frac{4r}{r+1}\int_{\mathbb{R}^3}\left(\left|\nabla\Big(n^{\frac{r+1}{2}}\Big)\right|^2+\left|\nabla\Big(p^{\frac{r+1}{2}}\Big)\right|^2\right)dx\\
 \leq&\  C(r+1)\int_{\mathbb{R}^3}(n^{r+1}+p^{r+1})dx+\frac{C}{r},\label{inL}
 \end{split}
 \eeq
 where $C$ is independent of $r$.
Based on the differential inequality \eqref{inL}, we can argue as in \cite[Supplement Lemma 5.1]{FFM} to obtain the uniform $L^\infty$  bounds for $n(t)$ and $p(t)$. The proof is complete.
\end{proof}

\textbf{Proof of Theorem \ref{connn}.}  Lemma \ref{unieee} yields the uniform-in-time $L^r$ estimates \eqref{unie}. Then the conclusion of Theorem \ref{connn} follows from the same argument as in \cite[Theorem 4.1]{WMZ}. The proof is complete.


\bigskip
 \textbf{Acknowledgments.} The authors are
grateful to the referees for their very helpful comments and suggestions.
 H. Wu was partially supported by
NSF of China 11001058, SRFDP and NSF of Shanghai 10ZR1403800. J. Jiang was partially
supported by NSF of China 11201468.

\end{document}